\documentclass[a4paper,10pt]{article}
\usepackage{amssymb,amsmath,amsthm}
\usepackage{amsfonts}
\usepackage{mathrsfs} 
\usepackage{dsfont} 
\usepackage{bbold}  
\usepackage[english]{babel}
\usepackage{hyperref}  
\usepackage{tensor}

\usepackage{color}
\usepackage{esint}
\numberwithin{equation}{section}

\newcommand{\h}{\eta}  

\DeclareMathOperator{\meas}{meas}
\DeclareMathOperator{\diag}{diag} 
\DeclareMathOperator{\dist}{dist} 
\newcommand{\uno}{\mathds{1}}

\newtheorem{theorem}{Theorem}[section]
\newtheorem{proposition}[theorem]{Proposition}
\newtheorem{lemma}[theorem]{Lemma}
\newtheorem{corollary}[theorem]{Corollary}

\newtheorem{remark}[theorem]{Remark}
 
\newtheorem{remarks}[theorem]{Remark}
\newtheorem{definition}[theorem]{Definition}
\newcommand{\be}{\begin{equation}}
\newcommand{\ee}{\end{equation}}

\newcommand{\om}{\omega}

\newcommand{\e}{\varepsilon}

\newcommand{\al}{\alpha}
\newcommand{\de}{\delta}

\newcommand{\Norma}{\norma \norma}
\newcommand{\R}{\mathds R}

\newcommand{\Z}{\mathds Z}

\newcommand{\N}{\mathds N}
\newcommand{\T}{\mathds T}

\renewcommand{\b }{\beta }
\newcommand{\s }{\sigma }
\newcommand{\ii }{{\rm i} }

\newcommand{\g }{\gamma}

\newcommand{\m }{\mu }
\newcommand{\vphi}{\varphi }

\newcommand{\norma}{| \!\! |} 

\newcommand{\ph}{\varphi}

\def\ba{\begin{aligned}}
\def\ea{\end{aligned}}
\def\beginm{\begin{multline}}
\def\endm{\end{multline}}

\newcommand{\ZZZ}{\mathds{Z}} 
\newcommand{\CCC}{\mathds{C}} 
\newcommand{\NNN}{\mathds{N}} 
 
\newcommand{\RRR}{\mathds{R}} 
\newcommand{\TTT}{\mathds{T}}

\newcommand{\calA}{{\mathcal A}} 
\newcommand{\BB}{{\mathcal B}}

\newcommand{\CCCC}{{\mathcal C}}

\newcommand{\calF}{{\mathcal F}} 
\newcommand{\calG}{{\mathcal G}} 
 
\newcommand{\calI}{{\mathcal I}}

\newcommand{\LL}{{\mathcal L}} 
\newcommand{\MM}{{\mathcal M}} 
\newcommand{\NN}{{\mathcal N}}

\newcommand{\RR}{{\mathcal R}}

\newcommand{\gota}{{\mathfrak a}}

\newcommand{\gotd}{{\mathfrak d}} 
\newcommand{\gote}{{\mathfrak e}}

\newcommand{\gots}{{\mathfrak s}}

\newcommand{\gotG}{{\mathfrak G}}

\newcommand{\gotM}{{\mathfrak M}}

\newcommand{\pD}{\zeta}
 \newcommand{\pippa}{\mathtt c}
  
\newcommand{\ol}{\overline} 
 
\newcommand{\Fullbox}{{\rule{2.0mm}{2.0mm}}} 
 
\newcommand{\EP}{\hfill\Fullbox\vspace{0.2cm}} 
\newcommand{\prova}{\noindent{\it Proof. }} 
\newcommand{\io}{\infty} 
 
\newcommand{\x}{\xi} 

\newcommand{\ka}{\kappa} 
 
\newcommand{\f}{\varphi} 
 
\newcommand{\del}{\partial}

\newcommand{\av}[1]{\langle #1 \rangle}

\newcommand{\oo}{{\omega}} 
\newcommand{\la}{{\lambda}}

\newcommand{\BBB}{\boldsymbol{B}} 
\newcommand{\bvert}{\boldsymbol{\vert}}

 \def\tilde#1{\widetilde{#1}}

\def\ins#1#2#3{\vbox to0pt{\kern-#2 \hbox{\kern#1 #3}\vss}\nointerlineskip}

\begin{document}

\title{{\bf Quasi-periodic solutions for the forced Kirchhoff equation on $\T^d$}}

\date{}

 \author{Livia Corsi$^{1}$, Riccardo Montalto$^{2}$
 \vspace{2mm} 
\\ \small 
 $^{1}$ School of Mathematics, Georgia Institute of Technology, 686 Cherry St. NW, Atlanta GA, 30332, USA
\\ \small 
$^2$  Institut f\"ur Mathematik, Universit\"at Z\"urich
Winterthurerstrasse 190, CH-8057 Z\"urich, CH
\\ \small 
E-mail: lcorsi6@math.gatech.edu, riccardo.montalto@math.uzh.ch}

\maketitle

\noindent
{\bf Abstract:}
In this paper we prove the existence of small-amplitude quasi-periodic solutions with Sobolev regularity, for the $d$-dimensional forced Kirchhoff equation
with periodic boundary conditions. This is the first result of this type for a quasi-linear equations in high dimension.
The proof is based on a Nash-Moser scheme in Sobolev class and a regularization procedure combined with a multiscale analysis in order to solve the linearized problem at any approximate solution. 
 \\[2mm]
\\[1mm]

\tableofcontents

\section{Introduction and main result}
In this paper we consider the forced Kirchhoff equation on the $d$-dimensional torus $\T^d$ 
\begin{equation}\label{main equation}
\partial_{tt} v - \Big( 1 + \int_{\T^d} |\nabla v|^2\, d x \Big) \Delta v = \delta f(\omega t, x)
\end{equation}
where $\delta> 0$ is a small parameter, $\omega := \lambda \bar \omega \in \R^\nu$, $\lambda \in {\cal I} := [1/2, 3/2]$, $\bar \omega$ a fixed diophantine vector, i.e. 
\begin{equation}\label{dio}
|\bar \omega \cdot \ell| \geq \frac{\gamma_0}{|\ell|^{\nu}}\,, \quad \forall \ell \in \Z^\nu \setminus \{ 0 \},
\end{equation}
and $f : \T^\nu \times \T^d \to \R$ is a sufficiently smooth function with zero average, i.e. 
\begin{equation}\label{zero mean f}
\int_{\T^{\nu + d}} f(\vphi, x)\, d \vphi\, d x = 0\,. 
\end{equation} 

Following \cite{B5,BB2,BCP} we assume also
\begin{equation}\label{diophquad}
\Big|\sum_{1\le i,j\le \nu}\ol{\om}_{i}\ol{\om}_{j}p_{ij} \Big|
\ge\frac{\g_{0}}{|p|^{\nu(\nu+1)}},
\qquad \forall\,p\in\ZZZ^{\nu(\nu+1)/2}\setminus\{0\}.
\end{equation}

Rescaling $v \mapsto \delta^{\frac13} v$, we see that \eqref{main equation} takes the form
 \begin{equation}\label{main equation eps}
\partial_{tt} v - \Big( 1 + \e \int_{\T^d} |\nabla v|^2\, d x \Big) \Delta v = \e  f(\omega t, x)\,, \quad \e:= \delta^{\frac23}\,. 
\end{equation}

Our aim is to prove the existence of quasi-periodic solutions of \eqref{main equation eps} for $\e$ small enough and $\lambda$ in a large subset of parameters in 
$\calI$. Since $\om$ is nonresonant, finding a quasi-periodic solution with frequency $\om$ is equivalent to find a torus embedding $\f\mapsto u(\f,\cdot)$
satisfying the equation $F(v) = 0$ where 
\begin{equation}\label{operator cal F}
{ F}(v) \equiv { F}(\lambda, v) := (\lambda \bar \omega \cdot \partial_\vphi)^2 v - \Big( 1 + \e \int_{\T^d} |\nabla v|^2\, d x \Big) \Delta v  -  \e  f(\vphi, x) 
\end{equation}
acting on the scale of real Sobolev spaces 
\begin{equation}\label{Sobolev real}
H^s = H^s(\T^{\nu + d} ) := \Big\{ v (\vphi, x) = \sum_{\begin{subarray}{c}
\ell \in \Z^\nu \\
j \in \Z^d
\end{subarray}} v_{\ell, j} e^{\ii \ell \cdot \vphi} e^{\ii j \cdot x}  \in L^2(\T^{\nu + d}) : \| v \|_s^2 :=  \sum_{\begin{subarray}{c}
\ell \in \Z^\nu \\
j \in \Z^d
\end{subarray}} \langle \ell, j \rangle^{2 s} |v_{\ell, j}|^2 < + \infty\Big\} 
\end{equation}
where $\langle \ell, j \rangle := {\rm max}\{1, |\ell|, |j| \}$.  Our main result is the following. 

\begin{theorem}\label{main}
There exist $q:=q(\nu,d)>0$, $s:=s(\nu,d)>0$ such that for any $f\in\CCCC^q(\TTT^\nu\times\TTT^d)$ satisfying \eqref{zero mean f} there exist 
$\e_0=\e_0(f,\nu,d)>0$ and for any $\e\in(0,\e_0)$ a Borel set $\CCCC_\e\subseteq \cal I$ with asymptotically full Lebesgue measure i.e.
$$
\lim_{\e\to0}\meas(\CCCC_\e)=1
$$
such that for any $\la\in\CCCC_\e$ there exists $u(\e,\la)\in H^s(\TTT^\nu\times\TTT^d)$ which is a zero for the functional $F$ appearing in \eqref{operator cal F}.
\end{theorem}

The Kirchhoff equation has been introduced for the first time in 1876 by Kirchhoff in dimension 1, without forcing term and with Dirichlet boundary conditions, to describe the transversal free vibrations of a clamped string in which the dependence of the tension on the deformation cannot be neglected. It is a quasi-linear PDE, namely the nonlinear part of the equation contains as many derivatives as the linear differential operator.

Concernig the existence of periodic solutions, Kirchhoff himself observed the existence of a sequence of {\it normal modes}, namely solutions of the form $v(t, x) = v_j(t) \sin(j x)$ where $v_j(t)$ is $2 \pi$-periodic. Under the presence of the forcing term $f(t, x)$ the {\it normal modes} do not persist\footnote{this is true except in the case where $f$ is uni-modal, i.e. $f(t, x) = f_k(t) \sin(k x)$ for some $k \geq 1$}, since, expanding $v(t, x) = \sum_{j} v_j(t) \sin(j x)$, $f(t, x) = \sum_{j} f_j(t) \sin(j x)$, all the components $v_j(t)$ are coupled.

The existence of periodic solutions for the forced Kirchhoff equation in any dimension has been proved by Baldi in \cite{Baldi Kirchhoff}, while the existence of quasi-periodic solutions in one space dimension under periodic boundary conditions has been proved in \cite{Montalto}.

Note that equation \eqref{main equation eps} is a quasi-linear PDE and it is well known that the existence of global solutions (even not periodic or quasi-periodic) for quasi-linear PDEs is not guaranteed, see for instance the non-existence results
in \cite{KM, Lax} for the equation $v_{tt} - a(v_x) v_{xx} = 0$, $a > 0$, $a(v) = v^p$, $p \geq 1$, near zero.

The existence of periodic solutions for wave-type equations with unbounded nonlinearities has been proved for instance in \cite{Rabinowitz-tesi-1967, B5, B2}.
For the water waves equations, which are fully nonlinear PDEs, we mention  \cite{Ioo-Plo-Tol, IP09, IP11, Alazard-Baldi}; see also \cite{Baldi-Benj-Ono} for fully non-linear Benjamin-Ono equations.  

The methods developed in the above mentioned papers do not work for proving the existence of quasi-periodic solutions.

The existence of quasi-periodic solutions for PDEs with unbounded nonlinearities has been developed by Kuksin \cite{K2} for KdV and then Kappeler-P\"oschel \cite{KaP}. This approach has been improved by Liu-Yuan \cite{LY0,LY} to deal with DNLS (Derivative Nonlinear Schr\"odinger) and Benjamin-Ono equations. These methods apply to dispersive PDEs  like KdV, DNLS but not to derivative wave equation (DNLW) which contains first order derivatives in the nonlinearity. KAM theory for DNLW equation has been recently developed by Berti-Biasco-Procesi in \cite{BBiP1, BBiP2}. Such results are obtained via a KAM-like scheme which is based on the so-called {\it second Melnikov conditions} and provides also the {\it linear stability} of the solutions.


The existence of quasi-periodic solutions can be also proved by imposing only {\it first order Melnikov conditions} and the so-called {\it multiscale approach}.
This method has been developed, for PDEs in higher space dimension, by Bourgain in \cite{Bo1,B3, B5} for analytic NLS and NLW,
extending the result of Craig-Wayne \cite{CW} for 1-dimensional wave equation with bounded nonlinearity. Later,
this approach has been improved by Berti-Bolle \cite{BB1, BB2} for NLW, NLS with differentiable nonlinearity and by Berti-Corsi-Procesi \cite{BCP} on compact Lie-groups. 

This method is especially convenient in higher space dimension since the second order Melnikov conditions are violated, due to the high multiplicity of the eigenvalues.
The drawback is that the linear stability is not guaranteed. Indeed there are very few results concerning the existence and linear stability of quasi-periodic solutions in the case of multiple eigenvalues. We mention \cite{ChierchiaYou, BKM} for the case of double eigenvalues and  \cite{EK, KAMbeam} in higher space dimension. 

 All the aforementioned results concern {\it semi-linear} PDEs, namely PDEs in which the order of the nonlinearity is strictly smaller than the order of the linear part.
 For quasi-linear (either fully nonlinear) PDEs, the first KAM results have been proved 
 by the {\it Italian team} in \cite{BBM-Airy, BBM-auto, BBM-mKdV, Giuliani, Feola-Procesi, Feola, Montalto, BertiMontalto, KAM ww f depth}.
 
 To the best of our knowledge all the results for quasi-linear and fully nonlinear PDEs are only in one space dimension. The result proved in this paper is the first one concerning the existence of quasi-periodic solutions for a quasi-linear PDE in higher space dimension.
 
 The reason why we achieve our result, whereas for other PDEs this is not possible (at least at the present time),
 is not merely technical and can be roughly explained as follows. 
 
 Almost all the literature about the existence of quasi-periodic solutions
 for dynamical systems in both finite and infinite dimension is ultimately related to a functional Newton scheme.
 It is well known that in the Newton scheme one has to solve the linearized problem, which in turn means that one has to invert
 the linearized functional. Such linearized functional is a liner operator acting on a scale of Hilbert spaces, hence one also
 needs appropriate bounds on the inverse in order to make the scheme convergent. Now, suppose that such linearized
 operator has the form ${\cal L} = \Delta + \e a(\vphi, x) \Delta$. In order to obtain bounds one wants to reduce this operator
 to constant coefficients up to a remainder (at least of order $0$). Passing to the Fourier side in space, the corresponding
 symbol is given by $ H(x, \xi) = |\xi|^2 + \e a(\vphi, x) |\xi|^2$ and hence reducing ${\cal L}$ to constant coefficients at leading order is equivalent
 to find a change of variables $(x, \xi) \mapsto (x', \xi')$ such that in the new variables the Hamiltonian $H(x, \xi)$  depends
 only on $\xi'$. In the one dimensional case this is always possible, whereas in dimension higher than one this is  possible only in very
 special cases, due to the Poincar\'e ``triviality'' Theorem stating that generically a quasi-integrable Hamiltonian is not integrable; see for instance \cite{Gal}.
 Of course there are some cases in which the Hamiltonian $H(x, \xi)$ is integrable
 (up to lower order terms); see for instance \cite{MontaltoIMRN,FGMP,albertino3}. Indeed in these cases
 the {\it complete reduction to constant coefficients} is achieved. However the three papers \cite{MontaltoIMRN,FGMP,albertino3}
 deal only with linear equations, whereas in the nonlinear case one has to fit the reducibility of the linearized operator
 with the Newton scheme. For instance, if in our  case one tries to follow the above scheme and reduce completely the linearized
 operator  (this is done in \cite{MontaltoIMRN}),
 one obtains a bound on the inverse of the linearized operator ${\cal L}(u)$ of the form
 $\|{\cal L}(u)^{- 1} h \|_s \lesssim_s \| h \|_{s + \sigma} + \| u \|_{2 s + \sigma} \| h \|_{s_0 + \sigma}$ for $s \geq s_0$, where
 $\sigma$ is a constant depending only on $\nu$ and $d$. It is well known  that a bound of this type is not enough for making 
 the Newton scheme convergent; see \cite{LZ}.
 
 In the present paper we overcome this difficulty as follows. First of all the highest order of our Hamiltonian symbol $H(x,\x)$ does not
 depend on $x$ so it is integrable; therefore we perform a reparametrization of time and we also apply a multiplication operator
 by a function depending only on time, and obtain a transformed operator of the form
 $$
 (\om\cdot \del_\theta)-\mu\Delta + \RR_2,
 $$
 where $\mu$ is a constant $\e$-close to $1$ and $\RR_2$ is a bounded operator satisfying decay bounds; see \eqref{fixed mu} and \eqref{R2}.
 Then we do not attempt a reduction scheme for the lower order term $\RR_2$ but rather use the multiscale approach.  
A priori this implies that we may not have informations about the linear stability of the solution we find; however the linear stability is 
obtained a-posteriori, namely here we prove the existence, then by linearizing on the found solution one can apply Theorem 1.2 of \cite{MontaltoIMRN}.
An a-posteriori approach of this type has been used for instance in \cite{CHP} for the NLS on $SU(2)$, $SO(3)$.

Out of curiosity we finally note that our remainder $\RR_2$ has a loss of regularity $\s$ which is due to change of variables
needed for the reduction up to order zero; see \eqref{R2}.
We find it interesting that a similar loss of reguarity appears for semi-linear PDEs when the space variable lives on a
compact Lie group instead of a torus; see (2.24c) in \cite{BCP} where such loss is denoted by $\nu_0$. 

The paper is organized as follows. After reducing the problem to the zero mean value functions, we introduce the scale of Hilbert spaces
and recall some of their properties. In Section \ref{brodino} we discuss some properties of the linearized operator $\LL(u)$, and we reduce it to constant
coefficients up to a remainder of order zero. We then discuss a Nash-Moser scheme converging on a set $A_\io$ defined in terms of the reduced
operator, and which in principle might be empty.
Afterwards in Section \ref{caffe} we introduce a subset $\CCCC_\io\subseteq A_\io$ where the multiscale approach can be used.
Finally we provide measure esitmates on another subest $\CCCC_\e\subseteq\CCCC_\io$, defined in terms of the final solution only.

\medskip

\noindent
{\bf Acknowledgements}. We warmly thank M. Procesi for carefully reading the manuscript, and for her comments and suggestions. 
L.C. was supported by NSF grant DMS-1500943. R.M. was supported by Swiss National Science Foundation, grant {\it Hamiltonian systems of infinite dimension}, project number: 200020--165537.


\section{Reduction on the zero mean value functions}

We follow \cite{Montalto}. Defining the projectors $\Pi_0, \Pi_0^\bot$ as the orthogonal projections 
$$
\Pi_0 v := v_0(\vphi) =  \frac{1}{(2 \pi)^d} \int_{\T^d} v(\vphi, x)\, d x \,, \quad \Pi_0^\bot := {\rm Id} - \Pi_0\,,
$$
and writing $v = v_0 + u$, $u := \Pi_0^\bot v$, $f = f_0 + g$, $g := \Pi_0^\bot f$, the equation ${F}(v) = 0$ (see \eqref{operator cal F}) is equivalent to 
\begin{equation}\label{sistem Pi0 Pi bot}
\begin{cases}
(\lambda \bar \omega \cdot \partial_\vphi)^2 u - \Big(1 + \e \int_{\T^d} |\nabla u|^2\, d x \Big) \Delta u - \e g = 0\,,\\
(\lambda \bar \omega \cdot \partial_\vphi)^2 v_0 - \e f_0 = 0\,. 
\end{cases}
\end{equation}

By \eqref{dio} and \eqref{zero mean f}, using that 
$$
\frac{1}{(2 \pi)^\nu}\int_{\T^\nu} f_0(\vphi)\, d \vphi = \frac{1}{(2 \pi)^{\nu + d}}\int_{\T^{\nu + d}} f(\vphi, x)\, d \vphi\, d x = 0
$$
the second equation in \eqref{sistem Pi0 Pi bot} is easily solved and we get
$$
v_0(\vphi) :=\e  (\lambda \bar \omega \cdot \partial_\vphi)^{- 2}  f_0\,. 
$$

Then we are reduced to look for zeroes of the nonlinear operator 
\begin{equation}\label{porca troia 0}
{\cal F}(u) \equiv {\cal F}(\lambda, u) := (\lambda \bar \omega \cdot \partial_\vphi)^2 u - \Big(1 + \e \int_{\T^d} |\nabla u|^2\, d x \Big) \Delta u - \e g
\end{equation}
acting on Sobolev spaces of functions with zero average in $x \in \T^d$, i.e. 
\begin{equation}\label{sobolev media nulla}
H^s_0 := \Big\{ u \in H^s : \int_{\T^d} u(\vphi, x)\, d x = 0  \Big\}\,. 
\end{equation}


\section{Function spaces, norms, linear operators}

Given a family of Sobolev functions $u(\vphi,  x; \lambda)$, $\lambda \in \Lambda \subset \R$,  we define the Sobolev norm $\norma \cdot \norma_s$ as 
\begin{equation}\label{Sobolev pesata}
\begin{aligned}
& \Norma u \Norma_s := \| u \|_s^{{\rm sup}} +  \| \partial_\lambda u \|_{s - 1}^{{\rm sup}}\,, \\
& \| u \|_s^{{\rm sup}} := \sup_{\lambda \in \Lambda} \| u (\cdot; \lambda)\|_s\,.
\end{aligned}
\end{equation}
If $\mu: \Lambda \to \R$, we define 
\begin{equation}\label{norma pesata scalare}
\Norma \mu \Norma := |\mu|^{\rm sup} + |\partial_\lambda \mu |^{\rm sup}\,, \quad |\mu|^{\rm sup} := \sup_{\lambda \in \Lambda} |\mu(\lambda)|\,. 
\end{equation}
Note that the classical interpolation result for $\norma \cdot \norma_s$ holds, i.e. given $u(\cdot; \lambda), v(\cdot; \lambda)$, $\lambda \in \Lambda$, one has 
\begin{equation}\label{inter}
\Norma u v \Norma_s \leq C(s) \Norma u \Norma_s \Norma v \Norma_{s_0} + C(s_0) \Norma u \Norma_{s_0} \Norma v \Norma_s\,, \quad s \geq s_0
\end{equation}
where we fix once and for all 
\begin{equation}\label{norma bassa fixed}
s_0 := \Big[ \frac{\nu + d}{2} \Big]  + 1
\end{equation}
and $[x]$ denotes the integer part of $x \in \R$. 

For any $N > 0$ let us define the spaces of trigonometric polynomials
\begin{equation}\label{trig}
E_N := {\rm span}\Big\{ e^{\ii (\ell \cdot \vphi + j \cdot x)} : 0 < |(\ell, j)| \leq N \Big\}\,
\end{equation}
and the orthogonal projector 
\begin{equation}\label{proiettore}
\Pi_N : L^2(\T^{\nu + d}) \to E_N \,, \quad \Pi_N^\bot := {\rm Id} - \Pi_N \,;
\end{equation}
of course the following standard smoothing estimates hold: 
\begin{equation}\label{stime smoothing}
\begin{aligned}
\Norma \Pi_N u \Norma_{s + \alpha} \leq N^\alpha \Norma u \Norma_s\,, \quad \Norma \Pi_N^\bot u \Norma_s \leq N^{- \alpha} \Norma u \Norma_{s + \alpha}\,. 
\end{aligned}
\end{equation}

Let us introduce the notations $\lesssim$ and $\lesssim_s$; we write
$a\lesssim b$ if there exists a constant $c=c(\nu,d,\g_0)$ such that $a < c b$, and $a\lesssim_s b$ if the constant depends also on $s$.

 We now recall some results concerning operators induced by diffeomorphism of the torus. 

\begin{lemma}{\bf }
\label{lemma:LS norms}
Let $\beta (\vphi; \lambda)$ satisfy $\Norma \beta \Norma_{s_0 + 1} \leq \delta$ for some $\delta$ small enough and $\omega = \lambda \bar \omega$ with $\lambda \in {\cal I}$.  
Then the composition operator 
$$
{\cal B} : u \mapsto {\cal B} u, \quad  
({\cal B}  u)(\ph,x) := u(\ph + \omega \beta(\vphi), x) \, , 
$$
satisfies 
\be\label{pr-comp1}
\| {\cal B} u \|_{s} \lesssim_s \| u \|_{s} 
+ \| \b \|_{s + s_0} \| u \|_1 \, , \qquad \mbox{ for all  } s \geq 1\,,
\ee
\begin{equation}\label{derivata lambda cal B}
\| (\partial_\lambda {\cal B}) u \|_s \lesssim_s \| u \|_{s + 1} + \Norma \beta\Norma_{s + s_0} \| u \|_2\,, \quad \forall s \geq 2\,. 
\end{equation}
Moreover the map $\vphi \mapsto \vphi + \omega \beta (\vphi)$ is invertible with inverse given by $\vartheta \mapsto \vartheta + \omega \breve \beta(\vartheta) $. The function $\breve \beta$ satisfies the estimate
\be\label{p1-diffeo-inv}
\Norma \breve \beta \Norma_{s} \lesssim_s  \Norma \b \Norma_{s + s_0} \, . 
\ee
\end{lemma}
\prova
The Lemma can be proved arguing as in the proof of Lemma B.4 in \cite{Baldi-Benj-Ono} (using also that by Sobolev embedding $\| \cdot \|_{{\cal C}^s} \lesssim \| \cdot \|_{s + s_0}$). The estimate on $\partial_\lambda {\cal B}$, follows by differentiating w.r. to $\lambda$, using the estimate \eqref{pr-comp1} and by applying the interpolation estimate \eqref{inter}. 
\EP

The following lemma follows directly by applying the classical Moser estimate for composition operators, see \cite{moser}.  
\begin{lemma}{\bf (Composition operator)} \label{Moser norme pesate}
Let $ f \in {\cal C}^{q}(\T^{\nu + d} \times B_K, \R )$, where $B_K := [- K, K]$ for some $K > 0$ large enough.  
If $u(\cdot; \lambda) \in H^s(\T^{\nu + d})$, $\lambda \in \Lambda$ is a family of Sobolev functions
satisfying $\| u \|_{s_0} \leq 1$. Then for any $s \geq s_0$ 
\begin{equation} \label{0811.10}
\Norma f(\cdot, u) \Norma_s \leq C(s,  f ) ( 1 + \Norma u \Norma_{s}) \, . 
\end{equation} 
\end{lemma}

\subsection{Linear operators on \texorpdfstring{$H^s_0$}{Hs} and matrices}
\label{sub.linearop}

Set $\ZZZ_*^d:= \ZZZ^d\setminus \{0\}$ and let $B,C\subseteq{\ZZZ^\nu\times\ZZZ^d_*}$. 
A bounded linear operator $ L : H_B^s \to H_C^s $ is represented, as usual,
 by a matrix in
\begin{equation}\label{lematrici}
\MM^B_C := 
\Big\{\big(M_{k}^{k'}\big)_{k\in C, k'\in B},\,
M_{k}^{k'}\in \CCC\Big\}.
\end{equation}

\begin{definition}\label{def:Ms} {\bf ($s$-decay norm)}
For any $M\in \MM^{B}_{C}$ 
we define its $s$-decay norm as
\begin{equation}\label{s-norm}
| M|_{s}^{2}:=\sum_{k\in{\mathfrak \ZZZ^{\nu+d}}}[M(k)]^{2}
\langle k\rangle^{2s}
\end{equation}
where, for $k=(\ell,j)$ $\av{k} := \max(1,|k| )=\max(1,|\ell|,|j|)$, 
\begin{equation}\label{maggiorante}
[M(k)]:=\left\{
\begin{aligned}
&\sup_{\substack{h-h'=k, h\in{C},\,h'\in{B}}} \big|M^{h'}_{h} \big|,
 &k\in {C}-{B}, \\
&\qquad 0, & k\notin {C}-{B}\,,
\end{aligned}
\right.
\end{equation}
If the matrix $M $ depends on a parameter $\lambda\in\Lambda\subseteq\RRR$, we define 
$$
\norma M\norma_s := |M|_s^{{\rm sup}} +  |\partial_\lambda M|_s^{\rm sup} \quad \text{where} \quad |M|_s^{{\rm sup}} := \sup_{\lambda \in \Lambda} | M(\lambda)|_s\,. 
$$
%
\end{definition}

\begin{remark}\label{estremo}
Note that if $M$ represent a multiplication operator by a function $a(\f,x)$ then
$$
| M|_s = \|a\|_s\qquad\mbox{and}\qquad \norma M\norma_s = \Norma a\Norma_s\,.
$$
\end{remark}

We have the following standard results; see for instance \cite{BB1} and references therein.

\begin{lemma}\label{lem.algebra.bvert} {\bf (Interpolation)}
For all $s\ge s_{0} $ there is $C(s)>1$ with $C(s_{0})=1$
such that, for any  subset $B,C,D\subseteq\ZZZ^\nu\times \ZZZ^d_*$
and for all $M_{1}\in\MM^{C}_{D}$, $M_{2}\in\MM^{B}_{C}$, one has
\begin{equation}\label{nanbnbna.bvert}
\bvert M_{1}M_{2}\bvert_{s}\le \frac{1}{2}\bvert M_{1}\bvert_{s_{0}}
\bvert M_{2}\bvert_{s}
+\frac{C(s)}{2}\bvert M_{1}\bvert_{s}\bvert M_{2}\bvert_{s_{0}}.
\end{equation}
In particular, one has the algebra property 
$ \bvert M_{1}M_{2}\bvert_{s}\le C(s) \bvert M_{1}\bvert_{s} \bvert M_{2}\bvert_{s} $. Similar estimates hold by replacing $| \cdot |_s$ with $\norma \cdot \norma_s$ if $M_1$ and $M_2$ depend on the parameter $\lambda$. 
\end{lemma}
Iterating the estimate of the above lemma one easily gets 
\begin{equation}\label{iterata}
| M^n |_s \leq C(s)^n | M |_s^{n - 1} | M |_{s_0}\,, \quad \forall n \in \N\,, \quad s \geq s_0\,. 
\end{equation}
If $M$ depends on the parameter $\lambda$, a similar estimate holds by replacing $| \cdot |_s$ with $\norma \cdot \norma_s$. 

\begin{lemma}\label{lem.questo}
For any $B,C\subseteq\ZZZ^\nu\times\ZZZ^d_*$,
let $M\in \mathcal M^B_C $.  Then 
\begin{equation}\label{soboh}
\| M h\|_s \leq C(s)\bvert M\bvert_{s_0}\|h\|_{s}+ C(s)\bvert M\bvert_{s}\|h\|_{s_0} \, , \quad \forall h\in H_B^s \, . 
\end{equation}
\end{lemma}


%
Of course all the results stated above hold replacing $|\cdot |_s$ by $\norma \cdot \norma_s$.


\section{The linearized operator}\label{brodino}

In this section we study the linearized operator ${\cal L}(u) := D_u {\cal F}(u)$ for any $u (\vphi, x; \lambda)$ which is ${\cal C}^\infty$ w.r.t.
 $(\vphi, x) \in \T^{\nu + d}$ and ${\cal C}^1$ w.r.t. the parameter $\lambda \in \cal I$.
 The linearized operator ${\cal L} : H^{s + 2}_0 \to H^s_0$, $s \geq 0$ has the form 
\begin{equation}\label{forma iniziale linearized}
\begin{aligned}
& {\cal L} = (\omega \cdot \partial_\vphi)^2 - \big(1 + a(\vphi) \big)  \Delta + {\cal R} \\
& a(\vphi) :=\e  \int_{\T^d} |\nabla u(\vphi, x)|^2\, d x \,,  \quad {\cal R}[h] := - 2 \Delta u \int_{\T^d} \Delta u \, h \, d x, \quad h \in L^2_0(\T^{\nu + d})\,. 
\end{aligned}
\end{equation}


\subsection{Reduction to constant coefficients up to the order zero}

In this section we prove the following Proposition.

\begin{proposition}\label{decay}
There exists  $\sigma = \sigma(\nu, d) > 0$  such that if 
\begin{equation}\label{ansatz u}
\Norma u \Norma_{s_0+\s } \leq 1\,,
\end{equation}
there exists $\delta \in (0, 1)$ such that 
if $\e \gamma_0^{- 1} \leq \delta$ then there exist two invertible changes of variables $\Phi_1, \Phi_2$ such that 
$$
\Phi_1 {\cal L} \Phi_2 = {\cal L}_2 = (\omega \cdot \partial_\vartheta)^2 - \mu \Delta + {\cal R}_2
$$ 
where $\mu$ is a constant and ${\cal R}_2$ is an operator of order $0$ satisfying the following properties.
The constant $\mu \equiv \mu (\lambda, u(\lambda))$ is ${\cal C}^1$ w.r.t. the parameter $\lambda$ and 
\begin{equation}\label{stima mu}
\Norma\mu - 1\Norma \lesssim \e \,, \quad |\partial_u \mu [h]|  \lesssim \e \| h \|_\sigma \,.
\end{equation}

The changes of variables $\Phi_1, \Phi_2$ are ${\cal C}^1$ w.r.t. the parameter $\lambda$ and they satisfy the tame estimates 
\begin{equation}\label{stime Phi 12 pm 1}
\begin{aligned}
& \| \Phi_1^{\pm 1} h \|_s , \| \Phi_2^{\pm 1} h\|_s \lesssim_s \| h \|_{s}
 + \| u \|_{s + \sigma} \| h \|_{s_0 } \,, \quad \forall s \geq s_0\,,  \\
& \| (\partial_\lambda\Phi_1^{\pm 1}) h \|_{s - 1} , \| (\partial_\lambda \Phi_2^{\pm 1}) h \|_{s - 1} \lesssim_s \| h \|_{s } + \Norma u \Norma_{s + \sigma}
\| h \| _{s_0 } \,, \quad \forall s \geq s_0.
\end{aligned}
\end{equation}

The remainder ${\cal R}_2$ is self-adjoint in $L^2$ and satisfies 
\begin{equation}\label{R2}
\begin{aligned}
& \norma{\cal R}_2\norma_s \lesssim_s \e (1 + \Norma u \Norma_{s + \sigma})\,, \quad \forall s \geq s_0\,, \\
& \norma\partial_u {\cal R}_2 [h]\norma_s \lesssim_s \e \Big( \Norma h \Norma_{s + \sigma} + \Norma u \Norma_{s + \sigma} \Norma h \Norma_{s_0 + \sigma} \Big)\,,
 \quad \forall s \geq s_0\,. 
\end{aligned}
\end{equation}
\end{proposition}


\subsubsection{Step 1: reduction of the highest order}

In this section we reduce to constant coefficients the highest order term $a(\vphi) \Delta$ in \eqref{forma iniziale linearized}. 
Given a diffeomorphism of the torus $\T^\nu \to \T^\nu$, $\vphi \mapsto \vphi + \omega \alpha (\vphi)$ we consider the induced operator 
\begin{equation}\label{definizione cal A}
{\cal A}h (\vphi, x) := h(\vphi + \omega \alpha(\vphi))
\end{equation} 
where $\alpha : \T^\nu \to \R$ is a small function to be determined. The inverse operator ${\cal A}^{- 1}$ has the form 
\begin{equation}\label{definizione cal A inv}
{\cal A}^{- 1} h(\vartheta, x) := h(\vartheta + \omega \breve \alpha(\vartheta), x)
\end{equation}
where $\vartheta \mapsto \vartheta + \omega \breve \alpha (\vartheta)$ is the inverse diffeomorphism of $\vphi \mapsto \vphi + \omega \alpha(\vphi)$. 
One has the following conjugation rules: 
\begin{equation}\label{regole coniugazione cal A A inverse}
\begin{aligned}
& {\cal A}^{- 1} a  {\cal A} = {\cal A}^{- 1}[a] \,, \quad  {\cal A}^{- 1} \circ \Delta \circ {\cal A} = \Delta \,, \\
& {\cal A}^{- 1}  (\omega \cdot \partial_\vphi) {\cal A} = {\cal A}^{- 1}\big[ 1 + \omega \cdot \partial_\vphi \alpha\big] \omega \cdot \partial_\vartheta \,, \\
& {\cal A}^{- 1}  (\omega \cdot \partial_\vphi)^2  {\cal A} = {\cal A}^{- 1}\big[ (1 + \omega \cdot \partial_\vphi \alpha)^2 \big] (\omega \cdot \partial_\vartheta)^2 
+ {\cal A}^{- 1}[(\omega \cdot \partial_\vphi )^2 \alpha] \omega \cdot \partial_\vartheta\,. 
\end{aligned}
\end{equation}

By \eqref{forma iniziale linearized}, \eqref{regole coniugazione cal A A inverse}, one has 
\begin{align}
{\cal A}^{- 1} {\cal L} {\cal A} & = {\cal A}^{- 1}\big[ (1 + \omega \cdot \partial_\vphi \alpha)^2 \big] (\omega \cdot \partial_\vartheta)^2 - 
{\cal A}^{- 1}[ 1 + a ] \Delta + {\cal A}^{- 1}[(\omega \cdot \partial_\vphi )^2 \alpha] \omega \cdot \partial_\vartheta + {\cal A}^{- 1} {\cal R} {\cal A} \,. 
\end{align}

We choose the function $\alpha$ so that the coefficient of $(\omega \cdot \partial_\vartheta)^2$ is proportional to the one of  the Laplacian $\Delta$, namely we want to solve 
\begin{equation}\label{bla bla}
(1 + \omega \cdot \partial_\vphi \alpha)^2 = \frac{1}{\mu} (1 + a) 
\end{equation}
for some constant $\mu \in \R$ to be fixed. Note that by \eqref{forma iniziale linearized}, \eqref{ansatz u}, one has that $a (\vphi) = O(\e)$, 
then for $\e $ small enough $\sqrt{1 + a}$ is well defined and of class $\CCCC^\io$. Then the equation \eqref{bla bla} can be written in the form 
\begin{equation}\label{blabla1}
\omega \cdot \partial_\vphi \alpha = \frac{1}{\sqrt{\mu}} \sqrt{1 + a} - 1\,. 
\end{equation}
and hence we choose $\mu$ so that the r.h.s. of \eqref{blabla1} has zero average, namely 
\begin{equation}\label{fixed mu}
\mu := \Big(\fint_{\T^\nu} \sqrt{1 + a(\vphi)}\, d \vphi \Big)^{2}\,.
\end{equation}

Now, using that $\omega = \lambda \bar \omega$ and $\bar \omega$ is diophantine, we choose 
\begin{equation}\label{definizione alpha}
\alpha := (\omega \cdot \partial_\vphi)^{- 1}\big[  \frac{1}{\sqrt{\mu}}  \sqrt{1 + a} - 1 \big]\,,
\end{equation}
and in this way, we obtain  
\begin{equation}\label{cal L1}
\begin{aligned}
& {\cal A}^{- 1} {\cal L} {\cal A} = \rho {\cal L}_1 \,, \quad \rho := {\cal A}^{- 1}[(1 + \omega \cdot \partial_\vphi \alpha)^2]\,, \\
& {\cal L}_1 := (\omega \cdot \partial_\vartheta )^2 - \mu \Delta + a_1 \omega \cdot \partial_\vartheta + {\cal R}_1 \,, \\
& a_1 := \rho^{- 1} {\cal A}^{- 1}[(\omega \cdot \partial_\vphi )^2 \alpha]  \,, 
& {\cal R}_1 := \rho^{- 1}{\cal A}^{- 1} {\cal R} {\cal A}\,. 
\end{aligned}
\end{equation}

\begin{lemma}
One has $\int_{\T^\nu} a_1(\vartheta)\, d \vartheta = 0$. 
\end{lemma}

\prova
By \eqref{cal L1} 
$$
a_1(\vartheta) = {\cal A}^{- 1}\Big[ \frac{(\omega \cdot \partial_\vphi )^2 \alpha}{(1 + \omega \cdot \partial_\vphi \alpha)^2} \Big] (\vartheta) 
= \frac{(\omega \cdot \partial_\vphi )^2 \alpha (\vartheta + \omega  \breve \alpha(\vartheta))}{(1 + \omega \cdot \partial_\vphi \alpha(\vartheta 
+ \omega  \breve \alpha(\vartheta)))^2}\,.
$$

Considering the change of variables $\vphi = \vartheta + \omega \breve \alpha(\vartheta)$, one gets 
\begin{align}
\int_{\T^\nu} a_1(\vartheta)\, d \vartheta & = \int_{\T^\nu} \frac{(\omega \cdot \partial_\vphi )^2 
\alpha(\vphi)}{(1 + \omega \cdot \partial_\vphi \alpha(\vphi))^2} (1 + \omega \cdot \partial_\vphi \alpha(\vphi))\, d \vphi \nonumber\\
& = \int_{\T^\nu} \frac{(\omega \cdot \partial_\vphi )^2 \alpha(\vphi)}{1 + \omega \cdot \partial_\vphi \alpha(\vphi)} \, d \vphi = 
\int_{\T^\nu} \omega \cdot \partial_\vphi \log\big( 1 + \omega \cdot \partial_\vphi \alpha(\vphi) \big)\, d \vphi = 0\,. 
\end{align}
\EP


\subsubsection{Step 2: reduction of the first order term}

The aim of this section is to eliminate the term $a_1(\vartheta) \omega \cdot \partial_\vartheta$ in the operator ${\cal L}_1 $ defined 
in \eqref{cal L1}. We conjugate ${\cal L}_1$ by means of a multiplication operator 
$$
{\cal B } : h \mapsto b(\vartheta) h
$$
where $b : \T^\nu \to \R$ is a function close to 1 to be determined, so that its inverse is given by 
$$
{\cal B}^{- 1} : h \mapsto b(\vartheta)^{- 1} h\,.
$$

One has the following conjugation rules: 
\begin{equation}\label{coniugio moltiplication}
\begin{aligned}
& {\cal B}^{- 1}\Delta {\cal B} = \Delta\,, \\
& {\cal B}^{- 1} \omega \cdot \partial_\vartheta {\cal B}  = \omega \cdot \partial_{\vartheta} + b(\vartheta)^{- 1} (\omega \cdot \partial_{\vartheta} b)\,, \\
& {\cal B}^{- 1} (\omega \cdot \partial_\vartheta)^2 {\cal B}  = (\omega \cdot \partial_{\vartheta})^2 +
 2 b(\vartheta)^{- 1} (\omega \cdot \partial_{\vartheta} b) \omega \cdot \partial_{\vartheta} + b(\vartheta)^{- 1} (\omega \cdot \partial_\vartheta )^2 b\,.
\end{aligned}
\end{equation}

By \eqref{cal L1}, \eqref{coniugio moltiplication} one gets 
\begin{align}
{\cal L}_2 & := {\cal B}^{- 1} {\cal L}_1 {\cal B} = (\omega \cdot \partial_\vartheta)^2 - \mu \Delta +
 \Big(  b(\vartheta)^{- 1} \omega \cdot \partial_{\vartheta} b + a_1(\vartheta)  \Big)  \omega \cdot \partial_\vartheta + {\cal R}_2 \label{primo cal L2}
\end{align}
where the remainder ${\cal R}_2$ is defined as 
\begin{equation}\label{def cal R2}
{\cal R}_2 := {\cal B}^{- 1} {\cal R}_1{\cal B} + b(\vartheta)^{- 1} (\omega \cdot \partial_\vartheta )^2 b + a_1(\vartheta) b(\vartheta)^{- 1} (\omega \cdot \partial_{\vartheta} b)\,. 
\end{equation}

In order to eliminate the term of order $\omega \cdot \partial_\vartheta$ one has to solve the equation 
\begin{equation}\label{equazione per b}
 b(\vartheta)^{- 1} \omega \cdot \partial_{\vartheta} b + a_1(\vartheta) = 0\,. 
\end{equation}
Since $ b(\vartheta)^{- 1} \omega \cdot \partial_{\vartheta} b = \omega \cdot \partial_\vartheta \log(b(\vartheta))$, the function $a_1$ has zero average,
 and recalling that $\omega = \lambda \bar \omega$ with $\bar \omega$ diophantine, the equation \eqref{equazione per b} can be solved by setting  
\begin{equation}\label{definizione b}
b(\vartheta) := {\rm exp}\Big(- (\omega \cdot \partial_\vartheta)^{- 1} a_1(\vartheta) \Big)\,. 
\end{equation}

Then ${\cal L}_2$ in \eqref{primo cal L2} has the final form 
\begin{equation}\label{L2}
{\cal L}_2 = {\cal  D} + {\cal R}_2\,, \quad {\cal D}={\cal D}(\la,u(\la)):= (\omega \cdot \partial_\vartheta)^2 - \mu \Delta \,,
\end{equation}
and the estimates \eqref{stima mu}-\eqref{R2}  follow similarly to \cite{Montalto}. Indeed they can be proved in an elementary way by using the explicit expressions for ${\cal R}_2, \Phi_1, \Phi_2, \mu$ found above and the estimate \eqref{inter}, Lemmata \ref{lemma:LS norms}, \ref{Moser norme pesate} and Remark \ref{estremo}.  

\begin{remark}\label{zero}
Note that for $u\equiv0$ one has $a=0$, $\mu=1$, $\al =1$, $\calA=\uno$, $\rho=1$, $a_1=1$, $b=1$, $\BB=\uno$ and hence
$$
\LL_2(0)=\LL(0) = (\om\cdot\del_{\vartheta})^2-\Delta\,.
$$
In particular $\RR_2(0)=0$. 
\end{remark}

\section{The Nash-Moser scheme.}\label{sec nash moser}

Here we prove the Nash-Moser scheme for parameters $\la$ in a set $A_\io$ (see below) which in principle might be empty; later we shall prove that
$A_\io$ contains the set $\CCCC_\e$ mentioned in Theorem \ref{main} and that $\CCCC_\e$ has asymptotically full measure.

For any $N > 0$ we decompose the operator ${\cal L} \equiv {\cal L}(u)$ as 
\begin{equation}\label{decomposizione cal L}
{\cal L}(u) = {\cal L}_N(u) + { \cal R}_N^\bot(u)
\end{equation}
where 
\begin{equation}\label{def cal DN cal RN}
\begin{aligned}
{\cal L}_N(u) & := \Phi_1(u)^{- 1} ( L_N (u)  + \Pi_N^\perp)\Phi_2(u)^{- 1}\,, \\
L_N(u) &:= { D}_N(\la,u(\la)) + { R}_N  \\
 { D}_N(\la,u(\la)) & := \Pi_N { \cal D}(\la,u(\la)) \Pi_N  \,, \\
 { R}_N(u) & := \Pi_N {\cal R}_2(u) \Pi_N \\
 { \cal R}_N^\bot(u) & := \Phi_1(u)^{- 1} \Pi_N^\bot{\cal L}_2(u) \Pi_N \Phi_2(u)^{- 1} + \Phi_1(u)^{- 1} \Pi_N{\cal L}_2(u) \Pi_N^\bot \Phi_2(u)^{- 1} \\
 & \quad  + \Phi_1(u)^{- 1} \Pi_N^\bot{\cal L}_2(u) \Pi_N^\bot \Phi_2(u)^{- 1} - \Phi_1(u)^{- 1} \Pi_N^\bot \Phi_2(u)^{- 1}\,.
\end{aligned}
\end{equation}

Note that, by applying the estimates \eqref{stime Phi 12 pm 1} and recalling \eqref{forma iniziale linearized}, the operator ${\cal R}_N^\bot$ satisfies 
\begin{equation}\label{stima cal R N modi alti}
\begin{aligned}
& \Norma {\cal R}_N^\bot h \Norma_{s_0} \lesssim N^{- \mathtt b} \big( \Norma h \Norma_{s_0 + \mathtt b + \sigma} + \Norma u \Norma_{s_0 + \mathtt b + \sigma}
 \Norma h \Norma_{s_0 + \sigma} \big)\,, \quad \forall \mathtt b > 0\,, \\
& \Norma {\cal R}_N^\bot h \Norma_{s} \lesssim_s  \Norma h \Norma_{s + \sigma} + \Norma u \Norma_{s  + \sigma} \Norma h \Norma_{s_0 + \sigma} \,, \quad \forall s \geq s_0\,. 
\end{aligned}
\end{equation}

Let $S>s_1>s_0+\s$ and consider $u \in {\cal C}^1({\cal I}, H^{s_1}_0)$ such that
\begin{equation}\label{us1}
\Norma u\Norma_{s_1} \leq 1\,;
\end{equation}
for any $\tau>0$, $\de\in(0,1/3)$ we define the set 
\begin{equation}\label{definizione cal GN}
\begin{aligned}
{\frak G}_N(u)  = {\frak G}_{N,\de,\tau}(u):= \Big\{ \lambda \in {\cal I}\; : \; &\forall\ s \in [s_1,S]\,,  \mbox{ one has} 
& | L_N(\lambda, u(\lambda))^{- 1} |_s \lesssim_s  
N^{ \frak a+\de (s - s_1)}(1+ \Norma u\Norma_{s+\s})\, 
 \Big\}\,,
\end{aligned}
\end{equation}
where $\frak a : =\tau+\de s_1$.

For any set $A \subset {\cal I}$ and $\eta > 0$ we define 
$$
{\cal N}(A, \eta) := \big\{ \lambda \in {\cal I} : {\rm dist}(\lambda, A) \leq \eta \big\}\,. 
$$
and let 
\begin{equation}\label{costanti nash-moser}
N_0 > 0\,, \quad N_n := N_0^{(3/2)^n}\,. 
\end{equation}

Let us introduce parameters $\ka_1$, $\ka_2$, $\ka_3$, satisfying
\begin{equation}\label{exponents}
\begin{aligned}
& \kappa_1 > \sigma \,, \quad \kappa_2 > {\rm max} \{ 3 \frak a + \frac{3}{2}(s_1 - s_0)+ 3 + \frac94 \kappa_1, 12 \frak a + 24 \} \,,  \\
& \kappa_3 > 6 \frak a + 6 + 3 \delta (S - s_1) + 3 \sigma + \frac32 \kappa_1 \,,\\
& (1 -  \delta) (S - s_1) > 2 \sigma + 2 + 2 \frak a + \frac23 \kappa_3 + \kappa_2\,. 
\end{aligned}
\end{equation}

Note one needs to impose the condition $0 < \delta < \frac13$ because the second and the third conditions are compatible only if 
$(1 - 3 \delta )(S - s_1) > 6 \frak a + 6 + \sigma + \kappa_1 $

\begin{theorem}\label{thm:nm1} {\bf (Nash-Moser)}
For $ \tau  $, $\de$, $ \ka_1$, $\ka_2$, $\ka_3$ ,$s_0 $, $ S >s_1>s_0+\s$, satisfying \eqref{exponents}, 
there are $ c $, $ \ol{N}_{0}$, such that, for all $N_0\ge \ol{N}_0$ and $ \e_0 $ small enough
such that 
\begin{equation}\label{piccoep} 
\e_0 N_0^S \leq c \, , 
\end{equation}
and, for all $\e\in[0,\e_0)$ a sequence
$\{u_n=u_{n}(\e,\cdot)\}_{n\ge0}\subset C^{1}(\mathcal I, H^{s_{1}}_0)$ such
that

\begin{itemize}

\item[(S1)$_{n}$] $u_{n}(\e,\la)\in E_{N_n}$, $u_{n}(0,\la)=0$, $\Norma u_n \Norma_{s_1} \leq 1$. 

\item[(S2)$_{n}$] For all $1\le i \le n$ one has
$\Norma u_{i}-u_{i-1}\Norma_{s_{1}}\le N_{i}^{- \kappa_1 }$. 

\item[(S3)$_{n}$] Set $u_{-1}:=0$ and define 
\begin{equation}\label{defAn}
A_{n}:=\bigcap_{i=0}^{n}
\gotG_{N_i}(u_{i-1}) \, . 
\end{equation} 
For $\la\in\NN(A_{n},N_n^{-\kappa_1/2})$
the function $u_{n}(\e,\la)$ satisfies $ \Norma {\cal F}(u_n) \Norma_{s_0} \leq C N_n^{- \kappa_2 } $.
\item[(S4)$_{n}$] For any $i = 1, \ldots, n$, $\Norma u_i \Norma_S \leq N_i^{\kappa_3}$. 

\end{itemize}

As a consequence, for all $\e\in[0,\e_0)$, 
 the sequence $\{u_{n}(\e,\cdot)\}_{n\ge0}$ converges uniformly in
$ C^{1}(\mathcal I, H^{s_1}_0)$ to $u_\e$ with
$ u_0(\la)\equiv 0$, at a superexponential rate
\begin{equation}\label{exponentialrate}
\Norma u_\e (\la)  - u_n (\la) \Norma_{s_1} \leq N_{n+1}^{-\kappa_1} \, , \quad \forall \la \in {\cal I} \, ,
\end{equation}
and for all $ \la\in A_{\io}:=\bigcap_{n\ge0}A_{n}$ one has $ \calF(\e,\la,u_\e (\la))=0 $. 
\end{theorem}

\subsection{Proof of Theorem \ref{thm:nm1}}

First of all we note that by differentiating the nonlinear operator ${\cal F}$ defined in \eqref{porca troia 0} by using \eqref{inter}, the following 
tame properties hold: for any $s\in[s_0,S]$ there
is $C=C(s)$ such that for any $u,h\in {\cal C}^1({\cal I}, H^{s}_0)$ with  $\Norma u\Norma_{s_0+2}\le1$ one has
\begin{itemize}

\item[(F1)] $\Norma\calF(\e,\la,u)\Norma_{s}\le C(s)(1+\Norma u\Norma_{s+2})$,

\item[(F2)] $\Norma D_{u}\calF(\e,\la,u)[h]\Norma_{s}\le C(s)(\Norma h\Norma_{s+2} + \Norma u\Norma_{s+2}\Norma h\Norma_{s_0+2})$,

\item[(F3)] $\Norma \calF(\e,\la,u+h)-F(\e,\la,u)-D_{u}F(\e,\la,u)[h]\Norma_{s}
\le C(s)(\Norma h\Norma_{s+2}\Norma h\Norma_{s_0+2} + \Norma u\Norma_{s+2}\Norma h\Norma_{s_0+2}^{2})$.

\end{itemize}

\begin{lemma}\label{lem.bbp2}
Let $\kappa > \mathfrak a + 2$ and $\Norma u\Norma_{s_{1}}\le 1$. For any  $\la\in {\cal N} \big(\gotG_N(u), 2 N^{- \kappa} \big)$, 
for $s \geq s_1$ there exists $\e_0 = \e_0(s)\in (0, 1)$ small enough such that if $\e \leq \e_0$,   the operator $L_N(\lambda, u(\lambda))$ is invertible and
\begin{subequations}
\begin{align}
& \norma L_{N}(u)^{-1} \norma_{s}  \lesssim_s N^{ 2 \frak a + 2 +  \delta(s - s_1) } (1 + \Norma u \Norma_{s + \sigma}) \,. 
\end{align}
\label{inv}
\end{subequations}
\end{lemma}

\prova
Let $\lambda \in \frak G_N(u)$ and $\lambda' \in {\cal I}$ so that $|\lambda - \lambda' | \leq 2 N^{- \kappa}$. We show by means of a Neumann series 
argument that $L_N(\lambda', u(\lambda'))$ is invertible, hence we want to bound $L_N(\e,\la',u(\la'))-L_N(\e,\la,u(\la))$.  By \eqref{stima mu} and \eqref{R2} we have
\begin{equation}\label{essi}
\begin{aligned}
| L_N(\e,\la',u(\la'))-L_N(\e,\la,u(\la)) |_{s}&\lesssim | \Pi_N ({\cal D}(\la,u(\la))-{\cal D}(\la',u(\la')))\Pi_N|_{s} \\
&\qquad+
| \Pi_N ({\cal R}_2(u(\la))-{\cal R}_2(u(\la')))\Pi_N|_{s} \\
&{\lesssim} (N^2 +  \e (1 + \Norma u \Norma_{s + \sigma})) |\la - \la'| \lesssim (N^2 +  \e (1 + \Norma u \Norma_{s + \sigma}))N^{- \kappa},
\end{aligned}
\end{equation}
so that for $s=s_0$, using that $s_0 + \sigma < s_1$ and $\Norma u \Norma_{s_1} \leq 1$ this reads
\begin{equation}\label{diffs0}
| L_N(\e,\la',u(\la'))-L_N(\e,\la,u(\la)) |_{s_0} \lesssim N^{- \kappa + 2}.
\end{equation}

Setting $A := L_N(\e,\la,u(\la))^{-1}(L_N(\e,\la',u(\la'))-L_N(\e,\la,u(\la)))$, by Neumann series one can write formally
$$
L_N(\lambda', u(\lambda'))^{- 1} = 
\sum_{n \geq 0} (- 1)^n A^n L_N(\lambda, u(\lambda))^{- 1}\,,
$$
and hence, using \eqref{essi}, \eqref{diffs0}, $\lambda \in \frak G_N(u)$
 and the interpolation estimate \eqref{nanbnbna.bvert}, we obtain 
\begin{equation}\label{stima A diff}
|A|_{s_0} \lesssim N^{2 + \frak a - \kappa} \,, \quad |A|_s \lesssim_s N^{\frak a + \delta (s - s_1) + 2 - \kappa} \big( 1 + \Norma u  \Norma_{s + \sigma} \big)\,,
\end{equation} 
so that by the estimate \eqref{iterata}, one obtains 
\begin{equation}\label{neu}
\begin{aligned}
|L_N(\e,\la',u(\la'))^{-1}|_s & \le \Big( \sum_{p\ge0} C(s)^p |A |_s |A|_{s_0}^{p - 1}\Big) |L_N(\e,\la,u(\la))^{-1}|_{s_0}  
+ \Big( \sum_{p\ge0} C(s_1)^p |A |_{s_0}^p\Big) |L_N(\e ,\la,u(\la))^{-1}|_{s} \\
& \lesssim_s N^{ \frak a + \delta(s - s_1)  } (1 + \Norma u \Norma_{s + \sigma})\,.  
\end{aligned}
\end{equation}

Now for any $\lambda \in {\cal N}\big( \frak G_N(u), N^{- \kappa} \big)$ by applying \eqref{def cal DN cal RN}, \eqref{stima mu}, \eqref{R2} one has 
\begin{equation}\label{stima partial lambda LN}
|\partial_\lambda L_N(\lambda, u(\lambda))|_s \lesssim_s N^2   + \Norma u \Norma_{s + \sigma} \,.
\end{equation} 

Finally, since $\partial_\lambda L_N(\lambda, u(\lambda))^{- 1} = - L_N(\lambda, u(\lambda))^{- 1} \partial_\lambda L_N(\lambda, u(\lambda))L_N(\lambda, u(\lambda))^{- 1}$, 
applying the estimates \eqref{neu}, \eqref{nanbnbna.bvert}, \eqref{stima partial lambda LN} one obtains that 
$$
|\partial_\lambda  L_N(\lambda, u(\lambda))^{- 1}|_s \lesssim_s N^{2 \frak a + 2 + \delta(s - s_1) } (1 + \Norma u \Norma_{s + \sigma})\,,
$$
so that the assertion follows.
\EP

The first step of the Nash-Moser algorithm is standard and uses the smallness condition \eqref{piccoep}.

Suppose inductively that $u_n$ is defined in such a way that the properties $(S1)_n-(S4)_n$ hold. We now define $u_{n + 1}$. We write 
\begin{align}
{\cal F}(u_n + h) & = {\cal F}(u_n) + D_u {\cal F}(u_n)[h] + {\cal Q}(u_n, h)
\end{align}
where 
\begin{equation}\label{def cal Q}
{\cal Q}(u_n, h) := {\cal F}(u_n + h) - {\cal F}(u_n) - D_u {\cal F}(u_n)[h]\,, 
\end{equation}
so that, using \eqref{decomposizione cal L} with $N = N_n$ and writing ${\cal F}(u_n) = \Pi_{N_{n + 1}} {\cal F}(u_n) + \Pi_{N_{n + 1}}^\bot {\cal F}(u_n)$ one gets 
\begin{align}
{\cal F}(u_n + h) & = {\cal F}(u_n) + {\cal L}_{N_{n + 1}}(u_n)[h]  + {\cal R}_{N_{n + 1}}^\bot (u_n)[h] + {\cal Q}(u_n, h)\,. \label{bla bla0}
\end{align}

Note that by applying Lemma \ref{lem.bbp2}, if $\lambda \in {\cal N}\big( A_{n + 1}, 2 N_{n + 1}^{- {\kappa_1}/{2}} \big)$ (recall \eqref{defAn}) 
the operator $L_{N_{n + 1}}(\lambda, u_n(\lambda)) : E_{N_{n + 1}} \to E_{N_{n + 1}}$ (recall \eqref{def cal DN cal RN}, \eqref{definizione cal GN}) is invertible, implying that $L_{N_{n + 1}}(\lambda, u_n(\lambda)) + \Pi_{N_{n + 1}}^\bot : H^s_0 \to H^s_0$ is invertible with $\norma \big( L_{N_{n + 1}}(\lambda, u_n(\lambda)) + \Pi_{N_{n + 1}}^\bot \big)^{- 1} \norma_s \leq \norma L_{N_{n + 1}}(\lambda, u_n(\lambda))^{- 1}\norma_s \lesssim_s N_{n + 1}^{ 2 \frak a + 2 +  \delta(s - s_1) } (1 + \Norma u_n \Norma_{s + \sigma})$. Since $\Phi_1(\lambda, u_n(\lambda))$ and $\Phi_2(\lambda, u_n(\lambda))$ are invertible for any $\lambda \in {\cal I}$ and satisfy the estimates 
\eqref{stime Phi 12 pm 1} then ${\cal L}_{N_{n + 1}}(\lambda, u_n(\lambda))$ is also invertible. By the estimates \eqref{stime Phi 12 pm 1}, 
the definition of the set $\frak G_{N_{n + 1}}(u_n)$, the estimate \eqref{soboh} and recalling that, by the inductive hypothesis $(S1)_n$
one has $\Norma u_n \Norma_{s_0 + \sigma} \leq \Norma u_n \Norma_{s_1} \leq 1$, 
we obtain 
\begin{equation}\label{accendino}
\Norma {\cal L}_{N_{n + 1}}(u_n)^{- 1}[h] \Norma_s \lesssim_s N_{n + 1}^{2 \frak a + 2} \Norma h \Norma_s + N_{n + 1}^{2 \frak a + 2 +  
\delta (s - s_1)} (1 + \Norma u_n \Norma_{s + \sigma}) \Norma h \Norma_{s_0}\,. 
\end{equation}

Let us now define, for $\lambda \in  {\cal N}\big(A_{n + 1}, 2 N_{n + 1}^{- {\kappa_1}/{2}} \big)$,  
\begin{equation}\label{def h n+1}
\widetilde h_{n + 1}(\lambda) := - \Pi_{N_{n + 1}} {\cal L}_{N_{n + 1}}(\lambda, u_n(\lambda))^{- 1} {\cal F}(\lambda, u_n(\lambda))\,, \quad \widetilde u_{n + 1} := u_n + \widetilde h_{n + 1}\,. 
\end{equation}

Plugging \eqref{def h n+1} into \eqref{bla bla0} one obtains  
\begin{equation}\label{F u n+1}
\begin{aligned}
{\cal F}(\widetilde u_{n + 1})  & = \Pi_{N_{n + 1}}^\bot {\cal F}(u_n) + {\cal R}_{N_{n + 1}}^\bot (u_n)[\widetilde h_{n + 1}] + {\cal Q}(u_n, \widetilde h_{n + 1})\,. 
\end{aligned}
\end{equation}

\medskip

\noindent
{\sc Estimate of $\widetilde h_{n + 1}$.} By applying \eqref{accendino}, using that $s_1 > s_0 + \sigma > s_0$, the property \eqref{stime smoothing} and $ \Norma u_n\Norma_{s_1} \leq 1$, one gets 
\begin{equation}\label{stima tilde h}
\begin{aligned}
\Norma \widetilde h_{n + 1} \Norma_{s_1} & \leq N_{n + 1}^{s_1 - s_0} \Norma {\cal L}_{N_{n + 1}}( u_n)^{- 1} {\cal F}(u_n) \Norma_{s_0}  \\
& \lesssim N_{n + 1}^{s_1 - s_0 + 2 \frak a + 2} \Norma {\cal F}(u_n) \Norma_{s_0}
 \stackrel{(S3)_n}{\lesssim} N_{n + 1}^{s_1 - s_0 + 2 \frak a + 2} N_n^{- \kappa_2}\,,  \\
\Norma \widetilde h_{n + 1} \Norma_S & \lesssim_S  N_{n + 1}^{2 \frak a + 2} \Norma {\cal F}(u_n) \Norma_S +
 N_{n + 1}^{2 \frak a + 2 +  \delta (S - s_1)} (1 + \Norma u_n \Norma_{S + \sigma}) \Norma {\cal F}(u_n) \Norma_{s_1} \\ 
& \stackrel{(F1), \eqref{stime smoothing}}{\lesssim_S } N_{n + 1}^{2 \frak a + 2 +  \delta (S - s_1) + \sigma} (1 + \Norma u_n\Norma_S)\,. 
\end{aligned}
\end{equation}

Let us consider a $C^\infty$ cut-off function $\psi_{n + 1}$ satisfying 
$$
\begin{aligned}
& {\rm supp}(\psi_{n + 1}) \subseteq {\cal N}\big( A_{n + 1}, 2 N_{n + 1}^{- \frac{\kappa_1}{2}} \big)\,, \quad 0 \leq \psi_{n + 1} \leq 1\,, \\
& \psi_{n + 1}(\lambda) = 1\,, \quad \forall \lambda \in {\cal N}\Big(A_{n + 1}, N_{n + 1}^{- \frac{\kappa_1}{2}}  \Big) \,.
\end{aligned}
$$
and define an extension of $\widetilde h_{n + 1}$ to the whole parameter space ${\cal I}$ as
$$
h_{n + 1}:= \psi_{n + 1} \widetilde h_{n + 1}\,, \quad u_{n + 1} := u_n + h_{n + 1}\,. 
$$

Using that $\Norma \psi_{n + 1} \Norma \lesssim N_{n + 1}^{\frac{\kappa_1}{2}}$ and by the estimates \eqref{stima tilde h} one has
\begin{subequations}\label{stime h}
\begin{align}
\Norma  h_{n + 1} \Norma_{s_1} & \lesssim N_{n + 1}^{s_1 - s_0 + 2 \frak a + 2 + \frac{\kappa_1}{2}} N_n^{- \kappa_2} \stackrel{\eqref{exponents}}{\lesssim} N_{n + 1}^{- \kappa_1}\,,
\label{acca1}  \\
\Norma  h_{n + 1} \Norma_S & \lesssim_S  N_{n + 1}^{2 \frak a + 2 +  \delta (S - s_1) + \sigma + \frac{\kappa_1}{2}} (1 + \Norma u_n\Norma_S)\,;
\label{accalta}
\end{align}
\end{subequations}
in particular $(S2)_{n+1}$ is satisfied.
Now
\begin{align}
\Norma u_{n + 1}\Norma_S & \lesssim_S \Norma u_n \Norma_S +  N_{n + 1}^{2 \frak a + 2 +  \delta (S - s_1) + \sigma + \frac{\kappa_1}{2}} (1 + \Norma u_n\Norma_S) 
\stackrel{(S4)_n}{\leq} C(S)N_{n + 1}^{2 \frak a + 2 + \delta(S - s_1) + \sigma + \frac{\kappa_1}{2}} N_{n }^{\kappa_3} \leq N_{n + 1}^{\kappa_3}
\end{align}
by \eqref{exponents} and by taking $N_0 = N_0(S) > 0$ large enough. Then also $(S4)_{n + 1}$ is proved. 

Now we estimate ${\cal F}(u_{n + 1})$ on the set ${\cal N}(A_{n + 1}, N_{n + 1}^{- \frac{\kappa_1}{2}})$. 
Using again that $\Norma u_n \Norma_{s_0+\s } \leq \norma u \norma_{s_1}<1$, one has  
\begin{align}
\Norma {\cal F}(u_{n + 1}) \Norma_{s_0} & \stackrel{\eqref{stime smoothing}, \eqref{stima cal R N modi alti}, (F3)}{\lesssim} 
N_{n + 1}^{- (S - s_0)} \Big( \Norma {\cal F}(u_n)\Norma_{S} + \Norma h_{n + 1} \Norma_{S + \sigma}  
+ \Norma u_n \Norma_{S + \sigma} \Norma h_{n + 1} \Norma_{s_1 }\Big)  + \Norma h_{n + 1} \Norma_{s_0 }^2  \nonumber\\
& \stackrel{(F1), \eqref{stime smoothing}, s_1 > s_0}{\lesssim} N_{n + 1}^{\sigma - (S - s_1)} \Big( 1 + \Norma u_n \Norma_S + \Norma h_{n + 1} \Norma_{S }  \Big)  
+ N_{n+1}^4 \Norma h_{n + 1} \Norma_{s_0}^2  \nonumber\\
& \stackrel{\eqref{stima tilde h}}{\lesssim} N_{n + 1}^{ 2 \sigma + 2  + 2 \frak a + (\delta - 1) (S - s_1) } \big(1 + \Norma u_n \Norma_S \big) 
+  N_{n + 1}^{4 \frak a +8} \Norma {\cal F}(u_n) \Norma_{s_0}^2 \nonumber\\
& \stackrel{(S3)_n, (S4)_n}{\lesssim} N_{n + 1}^{ 2 \sigma + 2  + 2 \frak a + (\delta - 1) (S - s_1) } N_n^{k_3} 
+ N_{n + 1}^{4 \frak a + 8} N_n^{- 2 \kappa_2} \leq N_{n + 1}^{- \kappa_2}  
\end{align}
by \eqref{exponents} and taking $N_0 = N_0(S) > 0$ large enough, hence proving $(S3)_{n+1}$.
Finally, 
by using a telescoping argument $u_{n + 1} = \sum_{i = 0}^{n + 1} h_i$, one has 
$$
\Norma u_{n + 1}\Norma_{s_1 } \stackrel{(S2)_n }{\leq} \sum_{i = 0}^{n + 1}  N_i^{- \kappa_1} \leq 1
$$
since by taking $N_0 > 0$ is large enough, thus providing $(S1)_{n+1}$.

Clearly the sequence $(u_n)_{n \in \N}$ is a Cauchy sequence in ${\cal C}^1({\cal I}, H^{s_1 }_0)$ and therefore the claimed statement follows. \EP

The proof of Theorem \ref{thm:nm1} is rather standard and follows the lines of the one in \cite{BBP,BCP}; however here we cannot apply directly
the aforementioned results because the subspaces $E_N$ in \eqref{trig} are not invariant under the change of variables $\calA$ appearing in \eqref{definizione cal A}.
We also mention that our truncation at the $n$-th step is not $N_0^{2^n}$ but rather $N_0^{\chi^n}$ with $\chi=3/2$; the reason for this choice is that, since
the subspaces $E_N$ are not invariant, we cannot apply the contraction Lemma at each step, but really the Newton scheme which converges only for $1<\chi<2$.

\section{Multiscale analysis}\label{caffe}

Our aim is to prove that the set $A_\io$ has asymptotically full measure; in order to do so, following \cite{BCP} we first prove that $A_\io$ contains another
set $\CCCC_\io$ and then we show that the set $\CCCC_\io$ contains another set $\CCCC_\e$ that has asymptotically full measure.

In order to do so,
in addition to the parameters $ 	\tau > 0 $, $ \delta \in (0, 1/3 ) $, $ \sigma $, $s_1$, $ s_0 $, $  S  $, $\ka_1,\ka_2,\ka_3$ satisfying
\eqref{exponents} needed in Theorem \ref{thm:nm1}, we now introduce other parameters 
 $\tau_{1} $,  $ \chi_0 $, $ \tau_0 $, $ C_1 $ and add the following constraints
\begin{equation}\label{esponenti}
\tau>\tau_0 \, ,\quad \tau_{1}>2\chi_0d \, , 
\quad
\tau>2\tau_{1}+d+\nu+1,\quad C_{1}\ge 2 \, , 
\end{equation}
then, setting $\ka:=\tau+d+\nu+s_{0}$, 
\begin{subequations}
\begin{align}
& \chi_0 (\tau-2\tau_{1}-d-\nu) > 3(\ka+(s_{0}+d+\nu)C_{1}),
\qquad \chi_0\de>C_{1},
\label{esponenti1.a} \\
&  s_{1} >3\ka+\s+2\chi_0(\tau_{1}+d+\nu)+C_{1}s_{0}. 
\label{esponenti1.b}
\end{align}
\label{esponenti1}
\end{subequations}
Note that no restrictions from above on $S'$ are required, i.e. it  could be $S'=+\io $.

Given $\Omega,\Omega'\subset\ZZZ^\nu\times \ZZZ^d_*$, we define
$$
{\rm diam}(\Omega):=\sup_{k,k'\in \Omega}\dist(k,k'),
\qquad
{\rm dist}(\Omega,\Omega'):=\inf_{\substack{k\in\Omega, k'\in\Omega'}}
|k-k'| \, , 
$$

\begin{definition}\label{regular} {\bf (Regular/singular sites)}
We say that the index $k=(\ell,j)\in\ZZZ^\nu\times \ZZZ^d_*$ is
\emph{regular} for a diagonal matrix  $ D $,  if $|D_{\ell,j} |\ge 1$, otherwise 
we say that $k$ is \emph{singular}.
\end{definition}

\begin{definition}\label{Ngood} {\bf ($N$-good/$ N$-bad matrices).}
Let $F\subset\ZZZ^\nu\times \ZZZ^d_*$ be such that ${\rm diam}(F)\le4N$
for some $N\in \NNN$.
We say that a matrix $A\in\MM_{F}^{F}$ is $N$\emph{-good}
if $A$ is invertible and for all $s\in[s_{0},s_{2}]$ one has
$$
\bvert A^{-1}\bvert_{s}\le N^{\tau+\de s}.
$$
Otherwise we say that $A$ is $N$\emph{-bad}.
\end{definition}

\begin{definition}\label{AN-reg} {\bf ($(A,N)$-regular, good, bad sites).}
For any finite  $E\subset\ZZZ^\nu\times \ZZZ^d_*$, let $ A = D + \e T \in \MM^E_E $ with 
 $ D := {\rm diag}( D_k ) $, $ D_k \in \CCC $. 
An index $k\in E$ is
\begin{itemize}

\item $(A,N)$-\emph{regular} if there
exists $F\subseteq E$ such that ${\rm diam}(F)\le4N$,
${\rm dist}(\{k\},E\setminus F)\ge N$ and the matrix $A^{F}_{F}$ is $N$-good.

\item $(A,N)$-\emph{good} if either it is regular for $D$ (Definition \ref{regular}) or it is
$(A,N)$-regular. Otherwise $k$ is $(A,N)$-\emph{bad}.
\end{itemize}
\end{definition}
The above definition could be extended to infinite $E$. 

Let $L$ be as in  \eqref{def cal DN cal RN}. 
Note that ${\cal D}$ in \eqref{L2} is represented by a diagonal matrix 
\begin{equation}\label{diag}
D(\lambda)  := {\rm diag}_{(\ell, j) \in\Z^\nu \times \Z^d_*} D_{\ell, j}(\lambda)\,, \quad D_{\ell, j}(\lambda) := - (\lambda \bar \omega \cdot \ell)^2 + \mu(\lambda) |j|^2\,. 
\end{equation}
Now for $\theta\in \RRR$ let us introduce the matrix
\begin{equation}\label{diagth}
D(\lambda,\theta)  := {\rm diag}_{(\ell, j) \in\Z^\nu \times \Z^d_*} D_{\ell, j}(\lambda,\theta)\,, \quad D_{\ell, j}(\lambda,\theta) :=
 - (\lambda \bar \omega \cdot \ell+\theta)^2 + \mu(\lambda) |j|^2\,,
\end{equation}
and denote
\begin{equation}\label{teta}
L(\e,\la,\theta,u):= D(\la,\theta)+\RR_2(u)\, .  
\end{equation}

\begin{lemma}\label{diniz}
For all $\tau>1, N>1, \lambda\in [1/2,3/2],\ell\in\ZZZ^\nu,j\in\ZZZ^d_*$ one has
\begin{equation}\label{inizio}
\{\theta\in\RRR\;:\: |D_{\ell, j}(\lambda,\theta)|\le N^{-\tau}\} \subseteq I_1\cup I_2\qquad
\mbox{ intervals with  }\meas(I_q)\le N^{-\tau}\,.
\end{equation}
\end{lemma}

\prova
A direct computation shows
$$
\{\theta \in \RRR\;:\; |D_{\ell,j}|\le N^{-\tau}_0\}=(\theta_{1,-},\theta_{1,+})\cup
(\theta_{2,-},\theta_{2,+})
$$
with
$$
\theta_{1,\pm}=\la\ol{\om}\cdot l +\sqrt{\mu|j|^2\pm N^{-\tau}}, 
\quad \theta_{2,\pm}=\la\ol{\om}\cdot l -\sqrt{\mu|j|^2\pm N^{-\tau}},
$$
and hence
$$
\meas((\theta_{q,-},\theta_{q,+})) =\frac{N^{-\tau} }{\sqrt{\mu|j|^2}}+ O(N^{-2\tau}), q=1,2.
$$
Note that by the estimate \eqref{stima mu}, $\mu \approx 1$ and $j \neq 0$ since we are working on the Sobolev space \eqref{sobolev media nulla}, so that the assertion follows.
\EP

For $ \tau_0 > 0 $, $N_0\ge 1$  we define the set
\begin{equation}\label{lambdabarra}
\ol{\mathcal I} := \ol{\mathcal I}(N_0,\tau_0):= \Big\{ \la\in\mathcal I\;:\; |  (\lambda \bar \omega \cdot \ell)^2 - |j|^2|
\ge N^{-\tau_{0}}_{0} \mbox{ for all }  k=(\ell,j)\in \ZZZ^\nu\times\ZZZ^d_*: \;|k|\le N_{0} \Big\}.
\end{equation}

In order to perform the multiscale analysis we need finite dimensional truncations
of such matrices. 
Given a parameter family of matrices $L(\theta)$ with $\theta\in\RRR$ and $N>1$ 
for any $k=(\ell,j)\in\ZZZ^\nu\times\ZZZ^d$ we denote by $L_{N,k}(\theta)$ (or equivalently $L_{N,\ell,j}(\theta)$) the 
sub-matrix of $L(\theta)$ centered at $ k $, i.e.
\begin{equation}\label{troncato}
L_{N,k}(\theta):=L(\theta)_{F}^F\,,\quad F:=\{k'\in \ZZZ^\nu\times\ZZZ^d_*\,: \; \dist(k,k')\leq N\}.
\end{equation}
If $\ell=0$, instead of the notation \eqref{troncato}  we shall use the notation
$$
L_{N,j}(\theta):=L_{N,0,j}(\theta) \, ,
$$
if also $j=0$ we write
$$
L_{N}(\theta):=L_{N,0}(\theta) ,
$$
and for $\theta=0$ we denote
$ L_{N,j}:=L_{N,j}(0) $. 

\begin{definition}\label{param.N-buoni}
{\bf ($N$-good/$ N$-bad parameters).} Let $\gote$ be large enough (to be computed).
We denote
\begin{equation}\label{bad}
B_{N}(j_0,\e,\la):= \Big\{\theta\in\RRR\,:\, L_{N,j_0}(\e,\la,\theta,u) \mbox{ is }N\mbox{--bad}\, \Big\}.
\end{equation}
A parameter $\la\in\mathcal I$ is $N$--good for $L $ if for any $j_0 \in \ZZZ^d$ one has
\begin{equation}\label{partizione}
B_{N}(j_0,\e,\la)\subseteq \bigcup_{q=1}^{N^{\gote}} I_{q} \, , \quad 
I_{q} \mbox{ intervals with } \meas(I_q) \le N^{-\tau_{1}} .
\end{equation}
Otherwise we say that $\la$ is
$N$--bad. We denote the set of $N$--good parameters as
\begin{equation}\label{param.good}
\calG_{N}=\calG_{N}(u):= \Big\{\la\in\mathcal I\,:\, \la\mbox{ is }N\mbox{--good for }L \Big\}.
\end{equation}
\end{definition}

The following assumption is needed for the multiscale Proposition \ref{multiscala}; we shall verify it later in Section \ref{separoL2}

\medskip

\noindent
{\bf Ansatz 1 (Separation of bad sites)} 
{\it 
There exist $ C_{1} >2$, 
$\hat{N}=\hat{N}(\tau_{0})
\in \NNN $ and $\hat {\mathcal I} \subseteq 
\overline{\mathcal I}$ (see \eqref{lambdabarra})
such that, for all $ N \ge \hat{N} $, and $ \| u \|_{s_1} < 1 $ (with $s_1$ satisfying \eqref{esponenti1.b}), if
$$
\la\in\calG_{N}(u)\cap \hat {\mathcal I},
$$
then for any $\theta\in\RRR$, for all $\chi\in[\chi_0,2\chi_0]$ and all $j_0\in\ZZZ^d$
the $(L,N)$-bad sites $k=(\ell,j)\in\ZZZ^\nu\times\ZZZ^d_*$
of $L=L_{N^{\chi},j_0}(\e,\la,\theta,u)$
admit a partition $\cup_{\b}\Omega_{\b}$ in disjoint
clusters satisfying
\begin{equation}\label{separazione}
{\rm diam}(\Omega_{\b})\le N^{C_{1}},
\qquad
{\rm dist}(\Omega_{\b_1},\Omega_{\beta_2})\ge N^{2},
\;\mbox{ for all }\b_1\ne\beta_2.
\end{equation}
}

%

\medskip

For $ N > 0 $, we  denote 
\begin{equation}\label{buoniautovalori}
\begin{aligned}
\calG_N^{0}(u):=\Big\{& \la\in\mathcal I\;:\;
\forall\; j_{0}\in \ZZZ^d  \mbox{ there is a covering } \\
B_{N}^{0}&(j_{0},\e,\la)\subset \bigcup_{q=1}^{N^{\gote}}I_{q}, \quad I_{q}=I_{q}(j_{0})\mbox{ 
intervals with }
\meas(I_{q})\le  N^{-\tau_1}
\Big\} 
\end{aligned}
\end{equation}
where
\begin{equation}\label{tetacattiviautovalori}
B_{N}^{0}(j_{0},\e,\la) := B_{N}^{0}(j_{0},\e,\la, u) := \Big\{\theta\in\RRR\;:\;
\|L_{N,j_{0}}^{-1}(\e,\la,\theta, u)\|_0>N^{\tau_1}\Big\} \, . 
\end{equation}
We also set
\begin{equation}\label{buoninormal2}
J_{N}(u):=\Big\{\la\in\mathcal I\;:\; \|L^{-1}_{N}(\e,\la,u)\|_{0}\le
N^{\tau_1} \Big\} \, .
\end{equation}
Under the smallness condition \eqref{piccoep}, Theorem  \ref{thm:nm1} applies, thus 
defining the sequence $ u_ n $ and the sets $ A_n $.
We  now introduce the sets
\begin{equation}\label{Cn}
\CCCC_{0}:=\hat{\mathcal I},\qquad
\CCCC_{n}:=\bigcap_{i=1}^n \calG^{0}_{N_i}(u_{i-1}) \bigcap_{i=1}^n J_{N_i}(u_{i-1})
 \cap \hat {\mathcal I}
\end{equation}
where $\hat {\mathcal I}$ is the one appearing in Proposition \ref{catene},  
$J_{N}(u)$  in \eqref{buoninormal2}, and  $\calG^{0}_{N}(u)$ in \eqref{buoniautovalori}. 

\begin{theorem}\label{thm:nm2}
Consider 
parameters satisfying  \eqref{exponents}, \eqref{esponenti}, \eqref{esponenti1}. 
Then
there exists $ \ol{N}_0 \in \NNN $, such that, for all $ N_0 \geq  \ol{N}_0 $ and $   \e \in [0, \e_0) $
with $ \e_0 $ satisfying \eqref{piccoep}, 
the following inclusions hold:
\begin{align*}
&(S5)_{0}\qquad\quad \;
\|u\|_{s_1}\le 1 \quad \Rightarrow\quad
\calG_{N_{0}}(u) ={\mathcal I} \\
&(S6)_0 \qquad\qquad \CCCC_{0}\subseteq A_{0},
\end{align*}
and for all $n\ge1$ (recall the definitions of $ A_n $ in \eqref{defAn}) 
\begin{align*}
&(S5)_{n}\qquad\quad \;
\|u-u_{n-1}\|_{s_1}\le N_{n}^{-\ka_1} \quad \Rightarrow\quad
\bigcap_{i=1}^{n}\calG^{0}_{N_i}(u_{i-1})\cap \hat {\mathcal I}\subseteq
\calG_{N_{n}}(u)\cap\hat {\mathcal I}, \\
&(S6)_{n}\qquad\qquad \CCCC_{n}\subseteq A_{n} \, .
\end{align*}
Hence  $ \CCCC_{\io}:=\bigcap_{n\ge0} \CCCC_{n} \subseteq A_\infty :=\bigcap_{n\ge0} A_n $. 
\end{theorem}

\subsection{Initialization}
\label{sec:inizio.nm2} 

Property  $(S5)_0$ follows from the following Lemma.

\begin{lemma}\label{iniz.1}
For all $\|u\|_{s_{1}}\le1$, $N\le N_0 $,   the set 
$\calG_{N}(u)=\mathcal I$.
\end{lemma}

\prova
We claim that, for any $\la\in[1/2,3/2]$ and any $j_{0}\in\ZZZ^d$, if (recalling the definition \eqref{diagth})
\begin{equation}\label{autov.grandi}
|D_{\ell,j}(\la,\theta)| > N^{-\tau_{1}},\quad\forall (\ell,j)\in \ZZZ^\nu\times \ZZZ^d_* \; \mbox{with}\; |(\ell,j-j_{0})|\le N \, , 
\end{equation}
then $ L_{N,j_{0}}(\e,\la,\theta)\;\mbox{ is }N\mbox{--good} $.
This implies that 
$$
B_{N}(j_{0},\e,\la)\subset \bigcup_{ |(l,j-j_{0})|\le N}\left\{ \theta\in\RRR\;:\; |D_{\ell,j}(\la,\theta)| \leq N^{-\tau_{1}}
\right\},
$$
which  in turn, by  Lemma \ref{diniz}, implies the thesis, see \eqref{partizione}, \eqref{param.good}, for some
$ \gote \geq d + \nu + 1 $.
The above claim  follows by a perturbative argument.
Indeed, recalling the definition \eqref{def cal DN cal RN}, for $\|u\|_{s_{1}} \le1 $, $ s_1 = s_2 + \s $, we
use \eqref{R2} to obtain
$$
\bvert (D_{N,j_{0}}^{-1}(\la,\theta))\bvert_{s_{2}}\bvert R_{N,j_{0}}(u)\bvert_{s_{2}} \leq \e C(s_1)
\bvert D_{N,j_{0}}^{-1}(\la,\theta) \bvert_{s_{2}}(1+\|u\|_{s_2 + \s})
\stackrel{\eqref{autov.grandi}}{\le}
\e N^{\tau_{1}}C(s_{1})\stackrel{\eqref{piccoep}}{\le} \frac{1}{2} \, . 
$$
Then we  invert $L_{N,j_{0}}$ by Neumann series and obtain
$$
\bvert L_{N,j_{0}}^{-1} (\e,\la,\theta)\bvert_{s}\le 2\bvert D_{N,j_{0}}^{-1}(\la,\theta) \bvert_{s}
\le 2 N^{\tau_{1}} {\le} N^{\tau+\de s}, \quad \forall s \in[s_{0},s_{2}]  \, , 
$$
by \eqref{esponenti}, which proves the claim.
\EP

\begin{lemma}\label{s60}
Property $(S6)_0$ holds.
\end{lemma}

\prova
Since $\hat{\mathcal I}\subset\ol{\mathcal I}$ it is sufficient to prove that $ \ol{\mathcal I} \subset A_0 $.
By the definition of $ A_0 $ in \eqref{defAn}, \eqref{definizione cal GN},  
we have to prove that
\begin{equation}\label{cio da dim}
\la\in\ol{\mathcal I} \quad \Longrightarrow \quad 
|L_{N_0}^{-1}(\e,\la,0)|_{s} \lesssim_s N_0^{\gota + \delta(s- s_1) }  \, , \ \forall s \in [s_1, S] \, .
\end{equation}
Indeed, if $\la\in\ol{\mathcal I}$ then $|D_{\ell,j}(\la)| \geq N_0^{-\tau_0}$, for all $|(\ell,j)|<N_0$, and so
$ \bvert D_{N_0}(\la)^{-1} \bvert_{s} \le  N_0^{\tau_0} $, $ \forall s $. 
Hence the assertion follows immediately by Remark \ref{zero} and \eqref{esponenti}.
\EP

\subsection{Inductive step} \label{sec:ter.nm2} 

By the Nash-Moser Theorem \ref{thm:nm1} we know that (S1)$_{n}$--(S4)$_{n}$ hold for all $n\ge0$.
Assume inductively that (S5)$_i $ and (S6)$_i $ hold for all $i \le n$. In order to prove
(S5)$_{n+1}$, we need the following {\em multiscale Proposition} \ref{multiscala} which
allows to deduce estimates on the $\bvert\cdot \bvert_s$--norm of the
inverse of $L$ from informations on the $L^2$-norm of the inverse $ L^{-1} $, 
the off-diagonal decay of $ L $, and separation properties of the bad sites.

\begin{proposition}\label{multiscala} {\bf (Multiscale)} Assume \eqref{esponenti}, \eqref{esponenti1}. For any $ \ol{s} > s_2 $, 
$\Upsilon>0$ there exists $\e_0=\e_0(\Upsilon,s_{2})>0$
and $N_{0}=N_{0}(\Upsilon, \ol{s})\in\NNN$ such that, for all $N\ge N_{0}$, $|\e|< \e_0$, $\chi\in [\chi_0,2\chi_0]$, 
$ E\subset \ZZZ^\nu\times \ZZZ^d_*$ with ${\rm diam}(E)\le 4N^{\chi}$, 
if the matrix $A=D+\e T\in\MM^{E}_{E}$ satisfies

\begin{itemize}
\item[{\rm (H1)}] $\bvert T\bvert_{s_{2}}\le \Upsilon$,
\item[{\rm (H2)}] $\|A^{-1}\|_{0}\le N^{\chi\tau_{1}}$,
\item[{\rm (H3)}] there is a partition $\{\Omega_{\beta}\}_{\beta}$
of   the $(A,N)$-bad sites (Definition \ref{AN-reg}) such that
$$
{\rm diam}(\Omega_{\beta})\le N^{C_{1}},
\qquad
{\rm dist}(\Omega_{\beta_1},\Omega_{\b_2})\ge N^{2},\;\mbox{ for }\b_1\ne\b_2,
$$
\end{itemize}

then the matrix $A$ is $N^{\chi}$-good and
\begin{equation}\label{stimabella}
\bvert A^{-1}\bvert_{s}\le\frac{1}{4}N^{\chi\tau}\left(N^{\chi\de s}+\e
\bvert T\bvert_{s}\right) \, , \quad \forall s \in [ s_0,\ol{s}] \, .  
\end{equation}
\end{proposition}
Note that the bound \eqref{stimabella} is much more than requiring that the matrix $A$ is 
$N^{\chi}$--good, since it holds also for $s > s_{2} $.

This Proposition is proved by ``resolvent type arguments'' and it coincides  essentially with \cite{BB1}-Proposition 4.1. The correspondences in the notations of this paper and \cite{BB1}
respectively are the following:
$(\tau,\tau_1,d+r,s_2, \ol{s})\rightsquigarrow (\tau',\tau,b,s_1, S)$, 
and, since we do not have a potential, we can fix  $\Theta=1$ in  Definition 4.2 of \cite{BB1}. 
Our conditions \eqref{esponenti}, \eqref{esponenti1} imply conditions (4.4) and (4.5) of \cite{BB1} 
for all $\chi\in [\chi_0,2\chi_0]$ and our (H1) implies the corresponding 
Hypothesis (H1) of \cite{BB1}  with $\Upsilon\rightsquigarrow 2\Upsilon$. The other hypotheses are the same.   
Although  the $s$--norm  in this paper is different, the proof of  \cite{BB1}-Proposition 4.1 relies only on abstract algebra and 
interpolation properties of the $s$--norm (which indeed hold also in this case -- see 
section \ref{sub.linearop}). Hence it can be repeated verbatim, 
full details can be found in arXiv:1311.6943.

\smallskip

Now,  we distinguish two cases:

\begin{itemize}

\item[{\bf case 1:}] $(3/2)^{n+1} \le \chi_0$. Then there exists $\chi\in[\chi_{0},2\chi_{0}]$ (independent of $n$) such that
\begin{equation}\label{case1}
N_{n+1}=\ol{N}^{\chi},\qquad \ol{N}:=[N_{n+1}^{1/\chi_0}]\in(N_{0}^{1/\chi},N_{0}) \, . 
\end{equation}
 This case may occur only in the first steps.

\item[{\bf case 2:}] $(3/2)^{n+1}>\chi_{0}$. Then there exists a unique $p\in[0,n]$ such that
\begin{equation}\label{case2}
N_{n+1}=N_{p}^{\chi},\qquad \chi=2^{n+1-p}\in[\chi_{0},2\chi_{0}) \, . 
\end{equation}
\end{itemize}

Let us start from {\bf case 1} for $n+1=1$; the other (finitely many) steps are identical.

\begin{lemma}\label{buonobarra}
Property (S5)$_{1}$ holds.
\end{lemma}

\prova
We have to prove that $ \calG_{N_1}^0(u_0) \cap \hat {\cal I} \subseteq {\cal G}_{N_1} (u) \cap \hat {\cal I} $
where $\|u-u_{0}\|_{s_{1}}\le N_{1}^{-\ka_1}$.  By Definition \ref{param.N-buoni} and \eqref{buoniautovalori}
it is sufficient to prove that, for all $j_{0}\in\ZZZ^d$, 
$$
B_{N_{1}}(j_{0},\e,\la,u)\subseteq B_{N_{1}}^{0}(j_{0},\e,\la,u_{0}),
$$
where we stress the dependence on $u,u_{0}$  in \eqref{bad}, \eqref{tetacattiviautovalori}. 
By the definitions \eqref{tetacattiviautovalori}, \eqref{bad} this amounts to prove that 
\begin{equation}\label{vweek}
\|L_{N_1,j_0}^{-1}(\e,\la,\theta,u_0)\|_{0}\le N_1^{\tau_1} \quad \Longrightarrow \quad  L_{N_1j_0}(\e,\la,\theta,u)\mbox{ is }N_1-\mbox{good} \, . 
\end{equation}
We first claim that $ \|L_{N_1,j_0}^{-1}(\e,\la,\theta,u_0)\|_{0}\le N_1^{\tau_1} $ implies
\begin{equation}\label{stimapotenteuno}
\bvert L^{-1}_{N_{1},j_{0}}(\e,\la,\theta,u_{0})\bvert_{s}\le
\frac{1}{4} N_{1}^{\tau}\left(N_{1}^{\de s}+
\bvert \RR_2(u_{0})\bvert_{s}\right) \stackrel{\eqref{R2}}{\le}
\frac{1}{4} N_{1}^{\tau}\left(N_{1}^{\de s}+
\e(1+\Norma u_0\Norma_{s+\s})\right) 
 \, , \quad \forall  s\in[s_0,S] \, .
\end{equation}
Indeed we may apply Proposition \ref{multiscala} to the matrix $A=L_{N_{1},j_{0}}(\e,\la,\theta,u_{0})$ with $\ol{s}=S$,
$ N = \ol{N} $, $ N_1 = \ol{N}^\chi $ and  $E=\{ |l|\le N_1,|j-j_0|\le N_1\}$. 
Hypothesis (H1)  
 follows by  \eqref{R2} and $ \| u_0 \|_{s_1 } \leq 1 $. Moreover (H2) is
$ \| L_{N_1,j_0}^{-1}(\e,\la,\theta,u_0)\|_{0}\le N_1^{\tau_1} $ . 
Finally (H3) is  implied by  Ansatz 1 
provided we take $N_0^{1/\chi_0}> \hat{N}( \tau_0)$ (recall \eqref{case1})
and noting that  $ \la \in \calG_{\overline N}(u_0) \cap \hat {\cal I} $
by Lemma \ref{iniz.1} (since $ \overline{N} \leq N_0 $ then $ \calG_{\overline N}(u_0) = \mathcal I  $). 
 Hence \eqref{stimabella} implies  \eqref{stimapotenteuno}. 
 
 We now prove \eqref{vweek}; we need to distinguish two cases.
 
 \medskip
 
 \noindent
 {\bf case 1. ($|j_0|>N_1^3$)}. We first show that $B_{N_1}^{0}(j_{0},\e,\la)\subset\RRR\setminus[-2N_1,2N_1]$. 
 Recall that if $A,A'$ are self-adjoint matrices, then their eigenvalues $ \mu_p (A) $, $ \mu_p (A') $ (ranked in nondecreasing order) 
satisfy
\begin{equation}\label{eq:7.10}
|\m_p (A) - \m_p (A') |\le \|A-A'\|_{0} \, .
\end{equation}

Threfore all the eigenvalues $\mu_{\ell,j}(\theta)$ of $L_{N_1,j_{0}}(\e,\la,\theta,u_0)$
are of the form 
\begin{equation}\label{diff auto}
\mu_{\ell,j}(\theta)
=\de_{\ell,j}(\theta)+O(\e\|\RR_2\|_{0}), 
\quad
\de_{\ell,j}(\theta) :=-(\oo\cdot \ell+\theta)^{2}+\mu(u_0)|j|^2 \, .
\end{equation}
Since $|\oo|_1=\la|\ol{\oo}|_1\le3/2$, $|j-j_{0}|\le N_1$, $|\ell|\le N_1$, one has 
$$
\de_{\ell,j}(\theta)\ge-\Big(\frac{3}{2}N_1+|\theta|\Big)^{2}+  N_1^2 > \frac{1}{2} N_1^2 \, \, , \quad \forall |\theta | < 2 N_1 \, . 
$$
and this implies $B_{N_1}^{0}(j_{0},\e,\la)\cap[-2N_1,2N_1]=\emptyset$. Hence the assumption
$\|L_{N_1,j_0}^{-1}(\e,\la,\theta,u_0)\|_{0}\le N_1^{\tau_1} $ implies $|\theta|<2N_1$. But then also the eigenvalues of 
$L_{N_1,j_0}(\e,\la,\theta,u)$ are big since they are also of the form
\begin{equation}\label{bah}
-(\oo\cdot \ell+\theta)^{2}+\mu(u)|j|^2+O(\e\|\RR_2\|_{0}). 
\end{equation}
But then this implies
$$
 L_{N_1j_0}(\e,\la,\theta,u)\mbox{ is }N_1-\mbox{good} \,.
$$

 \medskip
 
  \noindent
 {\bf case 2. ($|j_0|<N_1^3$)}.
Since  $\|u-u_{0}\|_{s_{1}}\le N_{1}^{-\ka_1}$ (recall that $\|u_0\|_{s_1}\le 1$ so $\|u\|_{s_1}\le 2$)
then 
\begin{equation}\label{perdio}
\begin{aligned}
\bvert L_{N_{1},j_{0}}(\e,\la,\theta,u_{0})&-L_{N_{1},j_{0}}(\e,\la,\theta,u)\bvert_{s_{2}}
\leq \bvert L_{N_{1},j_{0}}(\e,\la,\theta,u_{0})-L_{N_{1},j_{0}}(\e,\la,\theta,u)
\bvert_{s_{1}-\s}\\
&\le \bvert
(\mu(u_0)-\mu(u))\diag_{|j-j_0|,|\ell|<N_1}|j|^2 + R_N(u_0)-R_N(u)
\bvert_{s_1-\s}\\
&\lesssim N_1^6
\|u-u_0\|_{s_1}\le \frac{1}{2}
\end{aligned}
\end{equation}
By Neumann series  
and \eqref{stimapotenteuno} one has
$ \bvert L^{-1}_{N_{1},j_{0}}(\e,\la,\theta,u)\bvert_{s}\le 
N_{1}^{\tau+\de s} $
for all $s\in [s_{0},s_{2}]$, namely  $ L_{N_{1},j_{0}}(\e,\la,\theta,u)$ is $ N_1 $-good.
\EP

\begin{lemma}\label{buonobarraS6}
Property (S6)$_{1}$ holds.
\end{lemma}

\prova
Let $ \la \in \CCCC_1 := \calG^0_{N_{1}}(u_0)\cap J_{N_1}(u_0)  \cap \hat {\mathcal I} $, 
see \eqref{Cn}. By the definitions \eqref{defAn}, \eqref{definizione cal GN}, and (S6)$_0 $, 
in order to prove that $ \la \in A_1 $, 
it is sufficient to prove that $\la \in \gotG_{N_{1}}(u_0)$. Since 
$ \la \in 
J_{N_1}(u_0)  $ the matrix $ \| L_{N_1}^{-1}(\e,\la,u_0) \|_0 \leq N_1^{\tau_1} $ (see \eqref{buoninormal2}) and so
\eqref{stimapotenteuno} holds with $ j_0 = 0 $, $ \theta = 0  $. 
Hence $\la\in \gotG_{N_1}(u_0)$ 
\EP

Now we consider {\bf case 2}.

\begin{lemma}\label{lem:boh}
$
\bigcap_{i=1}^{n+1}\calG_{N_{i}}^{0}(u_{i-1})\cap \hat{\mathcal I} \subseteq \calG_{N_p}(u_n)\cap \hat{\mathcal I} $.
\end{lemma}

\prova
By $(S2)_{n} $ of Theorem \ref{thm:nm1} we get
$ \|u_{n}-u_{p-1}\|_{s_{1}}\le $ $ \sum_{i=p}^{n}\|u_{i}-u_{i-1}\|_{s_{1}} \le $ $ \sum_{i=p}^{n}N_{i}^{-\ka_1-1}{\le} $ $ N_{p}^{-\ka_1}\sum_{i=p}^{n}N_{i}^{-1}\le N_{p}^{-\ka_1} $.
 Hence  $ (S5)_{p} $ ($ p \le n $) 
implies 
$$
\bigcap_{i=1}^{n+1}\calG_{N_{i}}^{0}(u_{i-1})\cap \hat{\mathcal I}\subseteq
\bigcap_{i=1}^{p}\calG_{N_{i}}^{0}(u_{i-1})\cap \hat{\mathcal I}\stackrel{(S5)_{p}}{\subseteq} \calG_{N_{p}}(u_{n})\cap \hat{\mathcal I} 
$$
proving the lemma.
\EP

\begin{lemma}\label{induz:s5}
Property (S5)$_{n+1}$ holds.
\end{lemma}

\prova
Fix $\la\in \bigcap_{i=1}^{n+1}\calG_{N_{i}}^{0}(u_{i-1})\cap \hat{\mathcal I}                                         $.
Reasoning as in the proof of Lemma \ref{buonobarra}, it is sufficient to prove that, 
for all $ j_0 \in\ZZZ^d$, $\|u-u_{n}\|_{s_{1}}\le N_{n+1}^{-\ka_1}$, 
one has
\begin{equation}\label{cisiamo}
\|L^{-1}_{N_{n+1},j_{0}}(\e,\la,\theta,u_{n})\|_{0}\le N_{n+1}^{\tau_{1}} \quad \Longrightarrow \quad 
L_{N_{n+1},j_{0}}(\e,\la,\theta,u)\mbox{ is }N_{n+1}\mbox{--good} \, .
\end{equation}
We apply the multiscale Proposition \ref{multiscala} 
to the matrix $A=L_{N_{n+1},j_{0}}(\e,\la,\theta,u_{n})$ with $ N^\chi =  N_{n+1} $ and $ N = N_p $, see \eqref{case2}.
Assumption (H1) holds and (H2) is $ \|L^{-1}_{N_{n+1},j_{0}}(\e,\la,\theta,u_{n})\|_{0}\le N_{n+1}^{\tau_{1}} $. 
Lemma \ref{lem:boh}  implies that $ \la \in \calG_{N_p}(u_n)\cap \hat{\mathcal I} $ and therefore 
also (H3) is satisfied
by Ansatz 1. 
But then Proposition \ref{multiscala} implies
\begin{equation}\label{Ln+1ind}
\bvert L^{-1}_{N_{n+1},j_{0}}(\e,\la,\theta,u_{n})\bvert_{s}\le
\frac{1}{4} N_{n+1}^{\tau}\left(N_{n+1}^{\de s}+
\bvert \RR_2(u_{n})\bvert_{s}\right), \quad \forall s\in[s_{0},S] \, .
\end{equation}
Then we can follow word by word the proof of Lemma \ref{buonobarra}
(with $N_{n+1}$ instead of $N_1$, and $u_n$ instead of $u_0$), i.e. we
 separate the cases $|j_0|>N_{n+1}^3$ and $|j_0|\le N_{n+1}^3$
and the assertion follows.
\EP

\begin{lemma}\label{lem:finale}
Property (S6)$_{n+1}$ holds.
\end{lemma}

\prova
Again the proof follows
word by word the proof of Lemma \ref{buonobarraS6} 
with $N_{n+1}$ instead of $N_1$, and $u_n$ instead of $u_0$. 
\EP

Let us finally define the set
\begin{equation}\label{Cantorfinale}
 \CCCC_\e:= \bigcap_{n \geq 0}\bar \calG^0_{N_0^{2^n}} \cap\bar J_{N_0^{2^n}}\cap \tilde{\mathcal I} 
 \cap \ol{\mathcal I}
\end{equation}
where $  \tilde{\mathcal I}=\tilde{\mathcal I}(N_0)$
 is defined in Hypothesis 1, $ \ol{\mathcal I} $
in \eqref{lambdabarra}
and, 
for all $ N \in \NNN $, 
 \begin{equation}\label{buoninormal2finale}
\bar J_{N}:= \Big\{\la\in\mathcal I\;:\; \|L^{-1}_{N}(\e,\la,u_\e (\lambda))\|_{0}\le
 N^{\tau_1}/2 \Big\} \, , 
 \end{equation}
 \begin{equation}\label{buoniautovaloriinf}
 \begin{aligned}
 \bar \calG^0_{N}:=\Big\{& \la\in\mathcal I\;:\;
 \forall\; j_{0}\in \ZZZ^d  \mbox{ there is a covering } \\
 \bar B^0_{N}&(j_{0},\e,\la)\subset \bigcup_{q=1}^{N^{\gote}}I_{q},
 \mbox{ with }I_{q}=I_{q}(j_{0})\mbox{ intervals with }
 \meas(I_{q})\le  N^{-\tau_1}
 \Big\}
  \end{aligned}
 \end{equation}
  with 
  \begin{equation}\label{tetacattiviautovalorifinali}
\bar B^0_{N}(j_{0},\e,\la):= \Big\{\theta\in\RRR\;:\;
\|L_{N,j_{0}}^{-1}(\e,\la,\theta, u_\e(\lambda))\|_0>N^{\tau_1}/2\Big\} \, . 
 \end{equation}
 
 We have the following result.
 
 \begin{lemma}\label{CinftyCep}
 $ \CCCC_\e\subseteq \CCCC_\io $.
\end{lemma}

\prova
We claim that, for all $ n \geq 0 $,  the sets $\bar\calG^0_{N_n}\subseteq \calG^0_{N_n}(u_{n-1})$
and $ \bar J_{N_n} \subseteq J_{N_n}(u_{n-1}) $. 
These inclusions are a consequence of the super-exponential convergence \eqref{exponentialrate} of $ u_n $ to $ u_\e $. 
In view of the definitions \eqref{buoniautovaloriinf} and \eqref{buoniautovalori}, it is sufficient to prove that, $ \forall  j_0 $,
if  $\theta \notin \bar B_{N_n}^0(j_0,\e,\la)$ then 
$\|L_{N_n,j_0}^{-1}(\theta, u_{n-1})\|_0 \leq N_n^{\tau_1}$, namely $ \theta \notin B^0_{N_n}(j_0,\e,\la, u_{n-1}) $ (recall 
\eqref{tetacattiviautovalori}). Once again we have to distinguish two cases

\medskip

\noindent
 {\bf case 1. ($|j_0|>N_n^3$)}. In this case, arguing again as in the proof of Lemma \ref{buonobarra} one has
 $|\theta|<2N_n$, so the eigenvalues of $L_{N_n,j_0}(\theta, u_{n-1})$ are big and hence 
$\|L_{N_n,j_0}^{-1}(\theta, u_{n-1})\|_0 \leq N_n^{\tau_1}$.

\medskip

\noindent
 {\bf case 2. ($|j_0|\le N_n^3$)}.
One has  $ \|L_{N_n,j_0}^{-1}(\e, \la, \theta, u_{\e})\|_0 \leq N_{n}^{\tau_1} /2 $ by \eqref{tetacattiviautovalorifinali},
and so
 \begin{equation}\nonumber
\begin{aligned}
\|L_{N_n,j_0}^{-1}(\theta, u_{n-1})\|_0&\le \|L_{N_n,j_0}^{-1}(\theta, u_{\e})\|_0 \
\Big\| \Big(\uno+L_{N_n,j_0}^{-1}(\theta, u_{\e})(L_{N_n,j_0}(\theta, u_{n-1})- L_{N_n,j_0}(\theta, u_{\e}) )\Big)^{-1}\Big\|_0 \\
 &\le (N_{n}^{\tau_1} /2) \, 2 = N_{n}^{\tau_1}
\end{aligned}
\end{equation}
by Neumann series expansions. 
The inclusion $ \bar J_{N_n} \subseteq J_{N_n}(u_{n-1}) $
follow similarly. 
\EP

Theorem \ref{thm:nm2} and Lemma \ref{CinftyCep} are essentially Theorem 5.5 and Lemma 5.21 of \cite{BCP} respectively,
where \eqref{R2} implies Hypothesis 1 of \cite{BCP} with $\nu_0\rightsquigarrow \s$,
Lemma \ref{diniz} implies that Hypothesis 2 of \cite{BCP} is satisfied and Ansatz 1 here is the separation property of Hypothesys 4 in \cite{BCP}.
However we cannot directly apply the result of \cite{BCP} for the following reason.
The constant $\mu$ appearing in \eqref{diag} depends on the function at wich the linearized operator is computed; hence one has
$$
L_N(\e,\la,\theta,u) - L_N(\e,\la\theta,v) = (\mu(u)-\mu(v))\Delta + \RR_2(u) - \RR_2(v).
$$
The presence of the term $(\mu(u)-\mu(v))\Delta$ forces us to distinguish the cases $|j_0|$ large, where no small divisor appear, and 
$|j_0|$ small where one argues by Neumann series as in \cite{BCP}.

In what follows we are going to prove that Ansatz 1 is satisfied and later we shall
provide measure estimates for $\CCCC_\e$, thus concluding the proof of our main Theorem \ref{main}.

\section{Proof of Ansatz 1}\label{separoL2}

Given $\Sigma \subseteq \Z^\nu \times \Z^d_* $ we define for $\widetilde \jmath \in \Z^d_*$ the section  
$$
\Sigma^{(\widetilde{\jmath})}:=\{  k =(\ell ,\widetilde{\jmath})\in \Sigma\} \,.
$$

\begin{definition}\label{fibre}
Let $ \theta, \la $ be fixed and $ K > 1 $. We denote by $\Sigma_K$ any subset of singular sites of $D(\la,\theta)$
in $\ZZZ^\nu\times\ZZZ^d_*$ such that, for all $\widetilde{\jmath}\in \ZZZ^d_* $, the cardinality of the section $  \Sigma^{(\widetilde{\jmath})}_K $
satisfies  $  \# \Sigma^{(\widetilde{\jmath})}_K \le K$.
\end{definition}

\begin{definition}\label{gammachain} {\bf  ($ \Gamma$-Chain)} Let 
$\Gamma\ge 2 $. 
A sequence $k_{0},\ldots,k_{m}\in\ZZZ^\nu\times\ZZZ^d_*$ with
$k_{p}\ne k_{q}$ for $0\le p\ne q\le m$ such that
\begin{equation}\label{Gamma-chain}
\dist( k_{q+1},k_{q})\le \Gamma,\qquad
\mbox{ for all }q=0,\ldots,m-1,
\end{equation}
is called a $\Gamma$-chain of length $m$.
\end{definition}

\begin{proposition}\label{catene} {\bf (Separation of $\Gamma$-chains)}
There exists $ C=C(\nu,d)$ and, for any $ N_0 \ge 2 $
a set  $ \tilde{\mathcal I} = \tilde{\mathcal I}(N_0)$ defined as
\begin{equation}\label{lambdatilte}
\begin{aligned}
\tilde{\mathcal I} := \tilde{\mathcal I}(N_0)
:= &\Big\{\lambda\in[1/2,3/2]\;:\;   |P(\la\ol{\oo})|\ge\frac{N_0^{-1}}{1+|p|^{\nu(\nu+1)}}\; ,  \forall \, \mbox{non zero  polynomial} \\
&\; P(X)\in
\ZZZ[X_{1},\ldots,X_{\nu}]   \mbox{ of the form }
P(X)=p_{0}+\sum_{1\le i_{1}\le i_{2}\le \nu}p_{i_{1},i_{2}}X_{i_{1}}X_{i_{2}} \Big\} \, . 
\end{aligned}
\end{equation}
such that, 
for all $ \la\in \tilde{\mathcal I} $,  $ \theta \in \RRR $, and for all 
$ K, \Gamma $ with $ K \Gamma \geq N_0 $, any $\Gamma$-chain of
singular sites in $\Sigma_K$ as in Definition \ref{fibre},  has length
$ m\le (\Gamma K)^{C(\nu,d)} $.
\end{proposition}

\prova
The proof is a slight modification of  Lemma 4.2 of \cite{BB1} and Lemma 3.5 in \cite{BCP}. 
First of all, it is sufficient to bound the length of a
$\Gamma$-chain of singular sites for $D(\la,0)$. 
Then we consider the quadratic form 
\begin{equation}\label{quadratic-form}
Q:\RRR\times\RRR^{r}\to\RRR \, , \quad  Q(x,j):=-x^{2}+\mu|j|^{2},
\end{equation}
and the associated bilinear form  $ \Phi=- \Phi_{1}+\Phi_{2} $ where
\begin{equation}\label{split.b}
\Phi_{1}((x,j),(x',j')):=xx', \qquad
\Phi_{2}((x,j),(x',j')):=\mu j\cdot j' \, .
\end{equation}

For  a $\Gamma$-chain of sites $\{k_{q}=(\ell_{q},j_{q})\}_{q=0, \ldots, \ell}$ which are singular
for $D(\la,0)$ (Definition \ref{regular})
we have, recalling  \eqref{diag} and setting $x_{q} := \oo\cdot \ell_{q}$,   
$$
|Q(x_{q},j_{q})|<2, \qquad \forall q = 0,\ldots,\ell \, .
$$
Moreover, by \eqref{quadratic-form}, \eqref{Gamma-chain},  we derive
$|Q(x_{q}-x_{q_{0}},j_{q}-j_{q_0})|\le C|q-q_{0}|^{2}\Gamma^{2} $, $ \forall 0 \le q, q_0 \le m $, and so
\begin{equation}\label{stimameglio}
|\Phi((x_{q_{0}},j_{q_{0}}),(x_{q}-x_{q_{0}},j_{q}-j_{q_{0}}))|
\le C' |q-q_{0}|^{2}\Gamma^{2} \, . 
\end{equation}

Now 
we introduce the subspace of
$\RRR^{1+d}$ given by
$$
{\mathcal S}:={\rm Span}_{\RRR}\{(x_{q}-x_{q_{0}},j_{q}-j_{q_{0}})\;:
\;q=0,\ldots,m\}
$$
and denote by $\gots\le d+1$ the dimension of ${\mathcal S}$.
Let $\rho >  0 $ be  a small parameter specified later on. We distinguish two cases.

\noindent
{\bf Case 1.} {\it For all $q_{0}=0,\ldots,m$ one has
\begin{equation}\label{sottospazio.nlw-caso1}
{\rm Span}_{\RRR}\{(x_{q}-x_{q_{0}},j_{q}-j_{q_{0}})\;:
\;|q-q_{0}|\le \ell^{\rho},\;q=0,\ldots,m\}={\mathcal S}.
\end{equation}}
In such a case, we select a basis 
$ f_b :=(x_{q_{b}}-x_{q_{0}},j_{q_{b}}-j_{q_{0}})=(\oo\cdot\Delta \ell_{q_{b}},
\Delta j_{q_{b}})$, $b=1,\ldots,\gots$ of ${\mathcal S}$, where 
$\Delta k_{q_{b}}=(\Delta \ell_{q_{b}}, \Delta j_{q_{b}})$ satisfies
$\bvert \Delta k_{q_{b}}\bvert \le C\Gamma|q_{b}-q_{0}|\le C\Gamma m^{\rho}$. Hence
we have the bound
\begin{equation}\label{normabase}
\bvert f_{q_{b}}\bvert\le C\Gamma m^{\rho}, \qquad
b=1,\ldots,\gots.
\end{equation}
Introduce also the matrix $\Omega=(\Omega^{b'}_{b})_{b,b'=1}^{\gots}$ with 
$\Omega^{b'}_{b}:=\Phi(f_{b'},f_{b})$, 
that, according to \eqref{split.b}, we  write
\begin{equation}\label{chieomegone}
\Omega=\Bigl(- \Phi_{1}(f_{b'},f_{b})+\Phi_{2}(f_{b'},f_{b})
\Bigr)_{b,b'=1}^{\gots}=-X+Y,
\end{equation}
where $X^{b'}_{b} :=(\oo\cdot\Delta \ell_{q_{b'}})(\oo\cdot\Delta \ell_{q_{b}})$ and
$Y^{b'}_{b} :=\mu (\Delta j_{q_b'})\cdot (\Delta {j_{q_b}})$.
The matrix $Y$ has entries in $\mu \ZZZ$  and  the matrix
$X$ has rank $1$ since each  column is
\begin{equation}\nonumber
X^{b}=(\oo\cdot\Delta \ell_{q_{b}})
\begin{pmatrix}
\oo\cdot \Delta \ell_{q_{1}} \cr
\vdots \cr
\oo\cdot \Delta \ell_{q_{\gots}}\end{pmatrix}, \quad
b=1,\ldots,\gots.
\end{equation}
Then, 
since the determinant of a matrix with two collinear columns $X^{b},X^{b'}$, $b\ne b' $, is zero, we get
$$
\begin{aligned}
P(\oo):&=\mu^{d+1}{\rm det}(\Omega)=\mu^{d+1}{\rm det}(-X+Y)\\
&=\mu^{d+1}(\det(Y)-\det(X^{1},Y^{2},\ldots,Y^{\gots})-\ldots-
\det(Y^{1},\ldots,Y^{\gots-1},X^{\gots}))
\end{aligned}
$$
which is a quadratic polinomial as in \eqref{lambdatilte} with coefficients
$\le C(\Gamma m^{\rho})^{2(d+1)}$. Note that $P\not\equiv0$.
Indeed, if $P\equiv0$ then
\begin{equation}\nonumber
0=P(\ii\oo)=\mu^{d+1}\det(X+Y)=\mu^{d+1}\det(f_{b}\cdot f_{b'})_{b,b'=1,\ldots,\gots} \ne 0 
\end{equation}
because $\{f_{b}\}_{b=1}^{\gots}$ is a basis of ${\mathcal S}$. This contradiction proves that $ P \not\equiv 0 $. 
But then,  by \eqref{lambdatilte}, 
$$
\mu^{d+1}|\det(\Omega)|=|P(\oo)|\ge\frac{N_0^{-1}}{1+|p|^{\nu(\nu+1)}}\ge
\frac{N_0^{-1}}{(\Gamma m^{\rho})^{C(d,\nu)}} \, ,
$$
 the matrix $\Omega$ is invertible and 
\begin{equation}\label{inversaomegone}
|(\Omega^{-1})^{b'}_{b}|\le CN_0(\Gamma m^{\rho})^{C'(d,\nu)}.
\end{equation}
Now let 
$  {\mathcal S}^{\perp}:= {\mathcal S}^{\perp\Phi} := \{v\in\RRR^{s+1}\;:\;\Phi(v,f)=0,\;\forall\,f\in
{\mathcal S}\}$. Since $\Omega$ is invertible, the quadratic form
$\Phi_{{\mathcal S}}$ is non-degenerate and so
$\RRR^{d+1}={\mathcal S}\oplus{\mathcal S}^{\perp}$.
We denote  $\Pi_{{\mathcal S}}:\RRR^{d+1}\to{\mathcal S}$ the 
projector onto ${\mathcal S}$. Writing 
\begin{equation} \label{proj.esse}
\Pi_{{\mathcal S}}(x_{q_{0}},j_{q_{0}})=\sum_{b'=1}^{d+1}a_{b'}f_{b'} \, ,
\end{equation}
and since 
$f_{b}\in{\mathcal S}$, $ \forall b = 1,\ldots,\gots$, we get 
\begin{equation}\nonumber
w_b:=\Phi\big((x_{q_{0}},j_{q_{0}}),f_{b}\big)=
\sum_{b'=1}^{\gots}a_{b'}\Phi(f_{b'},f_{b})=\sum_{b'=1}^\gots\Omega_b^{b'} a_{b'}
\end{equation}
where $\Omega$ is defined in \eqref{chieomegone}.
The definition
of $f_{b}$, the bound \eqref{stimameglio} and \eqref{sottospazio.nlw-caso1} imply
$|w|\le C(\Gamma m^{\rho})^{2}$. Hence, by \eqref{inversaomegone},
we deduce $|a|=|\Omega^{-1}w|\le C' N_0(\Gamma m^{\rho})^{C(\nu,d)+2}$, whence, by
\eqref{proj.esse} and \eqref{normabase}, 
$$
|\Pi_{{\mathcal S}}(x_{q_{0}},j_{q_{0}})|\le N_0(\Gamma m^{\rho})^{C'(\nu,d)}.
$$
Therefore, for any $q_{1},q_{2}=0,\ldots,m $, one has
$$
|(x_{q_{1}},j_{q_{1}})-(x_{q_{2}},j_{q_{2}})|
=|\Pi_{{\mathcal S}}(x_{q_{1}},j_{q_{1}})-
\Pi_{{\mathcal S}}(x_{q_{2}},j_{q_{2}})|\le N_0(\Gamma m^{\rho})^{C_{1}(\nu,d)},
$$
which in turn implies $|j_{q_{1}}-j_{q_{2}}|\le N_0(\Gamma m^{\rho})^{C_{1}(r,d)}$
for all $q_{1},q_{2}=0,\ldots,m$. Since all the $j_{q}$ have $d$ components
(being elements of $\ZZZ^d_*$) they are at most
$CN_0^{d}(\Gamma m^{\rho})^{C_{1}(r,d)d}$. We are considering a $ \Gamma $-chain 
in $ \Sigma_K $ (see Definition \ref{fibre}) and so,
for each $q_{0} $, the number of $q \in \{0, \ldots, m \} $ such that $j_{q}=j_{q_{0}}$
is at most $K$ and hence
$$
m\le N_0^{d}(\Gamma m^{\rho})^{C_{2}(\nu,d)}K \leq (\Gamma K)^d  (\Gamma m^{\rho})^{C_{2}(\nu,d)}K \leq
 m^{\rho C_{2}(\nu,d)} (\Gamma K)^{d +C_{2}(\nu,d)} 
$$
because of the condition $\Gamma K \geq N_0 $,
Choosing $\rho<1/(2C_{2}(\nu,d))$ we  get $m \leq (\Gamma K)^{2(m +C_{2}(\nu,d))}$.

\noindent
{\bf Case 2.} There is $q_{0}=0,\ldots,m$ such that
\begin{equation}\nonumber
{\rm dim}(
{\rm Span}_{\RRR}\{(x_{q}-x_{q_{0}},j_{q}-j_{q_{0}})\;:
\;|q-q_{0}|\le m^{\rho},\;q=0,\ldots,m\})\le \gots-1.
\end{equation}
Then we  repeat the argument of Case 1 for the sub-chain
$\{(\ell_{q},j_{q})\;:\:|q-q_{0}|\le m^{\rho}\}$ and obtain a bound for $m^\rho $.
Since this procedure is applied at most $d+1$ times, at the end
we get a bound like $m\le(\Gamma K)^{C_{3}(\nu,d)}$.
\EP

\begin{corollary}\label{lama}
Ansatz 1 is satisfied.
\end{corollary}

The proof of Corollary \ref{lama} follows almost word by word Section 5.3 in \cite{BCP}.
However there is a minor issue to be discussed, namely that in Section 5.3 in \cite{BCP} it seems that one needs the index $j$ to be
in a lattice, whereas of course this is not the case in the present paper since we reduced to the zero mean valued functions.
However the lattice structure is needed only  in  Lemma 5.16 of \cite{BCP} (see Remark 5.17 of \cite{BCP}). In particular if we replace
Definition 5.14 of \cite{BCP} with  Definition \ref{stronglygood} below, the argument of \cite{BCP} can be repeated verbatim.

\begin{definition}\label{stronglygood}
A site $k=(\ell,j)\in\ZZZ^\nu\times\ZZZ^d  $ is

\begin{itemize}

\item $(L,N)$-strongly-regular if $L_{N,k}$ is $N$-good,
\item   $(L,N)$-weakly-singular if, otherwise, $L_{N,k}$ is $N$-bad,
\item $(L,N)$-strongly-good if either it is regular for $ D = D(\la, \theta ) $ (recall Definition \ref{regular}) or all the
sites $k'=(\ell',j')$ with
$\dist(k,k')\le N$ are $(L,N)$-strongly-regular. Otherwise $k$ is $(L,N)$-weakly-bad.

\end{itemize}
\end{definition}

\section{Measure estimates}

We conclude the proof of Thererm \ref{main} by showing that the set $\CCCC_\e$ has asymptotically full measure.

One proceeds differently for $|j_{0}| \ge 6N$ and $|j_{0}|< 6N$. 
We assume $N\ge N_0>0$ large enough and $\e \|\RR_2\|_0\le 1$.

\begin{lemma}\label{primapartizione} 
For all $j_{0}\in\ZZZ^d_*$, $| j_{0}|\ge 6N$, 
and for all $\la\in[1/2,3/2]$ one has
$$
\bar B_{N}^{0}(j_{0},\e,\la)\subset\bigcup_{q=1}^{N^{d+\nu+2}}I_{q}\, , 
\  \mbox { with } I_{q} = I_q(j_0)  \mbox{ intervals with } \meas(I_q) \le N^{-\tau_1} \, .
$$
\end{lemma}

\prova
First of all, as in the proof {\bf case 1} in Lemma \ref{buonobarra} we see that
 $\bar B_{N}^{0}(j_{0},\e,\la)\subset\RRR\setminus[-2N,2N]$. 
Now set
$ B_{N}^{0,+}:=\bar B_{N}^{0}(j_{0},\e,\la)\cap(2N,+\io) $,  $ B_{N}^{0,-}:=\bar B_{N}^{0}(j_{0},\e,\la)\cap(-\io, - 2 N) $.
Since
$$
\partial_{\theta}L_{N,j_{0}}(\e,\la,\theta)=
\diag_{{\substack|\ell|\le N,\\ |j-j_{0}|\le N}}-2(\oo\cdot \ell+\theta) \geq  N\uno,
$$
we apply Lemma 5.1 of \cite{BB2} 
with $\alpha= N^{-\tau_1}$, $\beta= N$ and $|E|\leq C N^{\nu+d}$ and obtain
\begin{equation}\nonumber
B_{N}^{0,-}\subset\bigcup_{q=1}^{N^{d+\nu+1}}I_{q}^{-}\ ,  \quad  I_{q}^{-}=I_{q}^{-}(j_{0}) \mbox{ intervals with }
\meas (I_q) \le  N^{-\tau_1} \, . 
\end{equation}
We can reason in the same way for $B_{N}^{0,+}$ and the lemma
follows.
\EP

Consider now $|j_{0}|<6N$. We  obtain a complexity
estimate for $\bar B_{N}^{0}(j_{0},\e,\la)$ by knowing the measure of the set
$$
\bar B^{0}_{2,N}(j_{0},\e,\la):=\left\{\theta\in\RRR\;:\;\|L_{N,j_{0}}^{-1}(\la,
\e,\theta)\|_{0}>N^{\tau_1}/2\right\}.
$$

\begin{lemma}\label{lem.controllobdoppio}
For all $|j_{0}|< 6 N$ and all $\la\in[1/2,3/2]$ one has
$$
\bar B^{0}_{2,N}(j_{0},\e,\la)\subset I_{N}:=[-  10\sqrt{d} N,  10\sqrt{d} N] 	 .
$$
\end{lemma}

\prova
If $|\theta|> 10\sqrt{d}  N$ one has
$|\oo\cdot \ell+\theta|\ge|\theta|-|\oo\cdot \ell|> (10\sqrt{\nu} - (3/2) )N > 8\sqrt{d} N $.
and then all the eigenvalues satisfy
$$
\mu_{\ell,j}(\theta)=-(\oo\cdot \ell+\theta)^{2}+\mu|j|^2+O(\e\|\RR_2\|_{0}) \le -62d N^{2} \, , \quad \forall |\theta| > 10\sqrt{d} N \, ,
$$
 proving the lemma.   
\EP

\begin{lemma}\label{secondocontrollo}
For all $|j_{0}|\le6N$
and all $\la\in[1/2,3/2]$ one has
$$
\bar B^{0}_{N}(j_{0},\e,\la)\subset\bigcup_{q=1}^{\hat{C}\gotM N^{\tau_1+1}}I_{q} \, , 
\ I_{q}=I_{q}(j_{0}) \mbox{ intervals with } \meas (I_q) \le N^{-\tau_1} 
$$
where $\gotM:=\meas(\bar B^{0}_{2,N}(j_{0},\e,\la))$ and $\hat{C}=\hat{C}(d)$.
\end{lemma}

\prova
This is Lemma 5.5 of \cite{BB2}, where our  exponent $\tau_1$ is denoted by $\tau$.
\EP

Lemmas \ref{lem.controllobdoppio} and \ref{secondocontrollo} imply that for all $\la
\in[1/2,3/2]$ the set $\bar B^{0}_{N}(j_{0},\e,\la)$ can be covered by 
$\sim N^{\tau_1+2}$ intervals of length $\le N^{-\tau_1}$. This estimate is not enough. 
Now we prove that for ``most" $\la$ the number of such intervals does not depend on $\tau_1$, by showing
that 
$\gotM = O( N^{\gote - \tau_1}) $ where $ \gote$ depends only on the dimensions (to be computed). To this purpose first  we  provide an estimate for the set
$$
\BBB^{0}_{2,N}(j_{0},\e):=\left\{(\la,\theta)\in[1/2,3/2]\times \RRR\;:\;
\|L_{N,j_{0}}^{-1}(\e,\la,\theta)\|_{0}>N^{\tau_1}/2\right\} \, . 
$$
Then in Lemma \ref{usofubini} we  use Fubini Theorem to obtain the desired bound for $\meas
(\bar B^{0}_{2,N}(j_{0},\e,\la))$.

\begin{lemma}\label{lem.misurabione.nlw}
For all $|j_{0}|< 6 N$ one has
$ {\rm meas}(\BBB^{0}_{2,N}(j_{0},\e))\lesssim  N^{-\tau_1+\nu+d+1}  $.
\end{lemma}

\prova
Let us introduce the variables
\begin{equation}\label{variabilinuove}
\pD=\frac{1}{\la^{2}},  \ \h=\frac{\theta}{\la}, \qquad
(\pD,\h)\in[4/9,4]\times [- 20\sqrt{d} N, 20\sqrt{d} N]=: [4/9,4] \times J_N,
\end{equation}
and set
$$
L(\pD,\h):= {\la^{-2}}L_{N,j_{0}}(\e,\la,\theta)=
\diag_{|\ell|\le N,|j-j_{0}|\le N}\Big(\big(-(\ol{\oo}\cdot \ell+\h)^{2} +
\pD\mu(\pD^{-1/2})|j|^2  \Big)+ \pD \RR_2(\e,1/\sqrt{\pD}).
$$
Note that, since $\norma \mu-1\norma\lesssim\e$, one has
\begin{equation}\label{minla}
\min_{j\in\ZZZ^d_*}\mu|j|^2\ge \frac{1}{2}.
\end{equation}
Then, except for $(\pD,\h)$ in a set of measure
$O(N^{-\tau_1+\nu+d+1})$ one has
\begin{equation}\label{stimotta}
\|L(\pD,\h)^{-1}\|_{0}\le N^{\tau_1}/8.
\end{equation}
Indeed
$$
\begin{aligned}
\partial_{\pD}L(\pD,\h)&=\diag_{|\ell|\le N,|j-j_{0}|\le N}\left(\mu(\pD^{-1/2})|j|^2 -\frac{1}{2}\pD^{-1/2} \del_\la\mu(\pD^{-1/2}))\right)
+  \RR_2(\e,1/\sqrt{\pD}) - \frac{1}{2} \pD^{-1/2}\partial_{\la}\RR_2
\stackrel{\eqref{minla}}  \ge \frac{1}{4},
\end{aligned}
$$
for $ \e $ small (we used that $\pD\in[4/9,4]$ and $|\del_\la\mu|<1/2$).
Therefore Lemma 5.1 of \cite{BB2} implies that for each $\h$,
the set of $\pD$ such that at least one eigenvalue of $L(\pD,\h)$ has
modulus $\le8 N^{-\tau_1}$, is contained in the union of $O(N^{d+\nu})$ intervals
with length $O(N^{-\tau_1})$ and hence has measure $\le O(N^{-\tau_1+d+\nu})$. Integrating
in $\h\in J_{N}$ we obtain \eqref{stimotta} except in a set with measure
$O(N^{-\tau_1+d+\nu+1})$. 
The same measure estimates hold in the original variables $ (\la, \theta) $ in \eqref{variabilinuove}. 
Finally \eqref{stimotta} implies
$$
\|L^{-1}_{N,j_{0}}(\e,\la,\theta)\|_{0}\le\la^{-2}N^{\tau_1}/8\le N^{\tau_1}/2,
$$
for all $(\la,\theta)\in[1/2,2/3]\times\RRR$ except in a set with
measure $\le O(N^{-\tau_1+d+\nu+1})$.
\EP

Note that the same argument can be used to show that 
\begin{equation}\label{cattivinormalta}
\meas([1/2,3/2]\setminus \bar \gotG_{N})\le N^{-\tau_1+d+\nu +1} 
\end{equation}
where $ \bar \gotG_N$ is defined in \eqref{buoninormal2finale}.

Define the set
\begin{equation}\label{FNj0}
\calF_{N}(j_{0}):=\left\{\la\in[1/2,3/2]\;:\;\meas(\bar B_{2,N}^{0}(j_{0},\e,\la))
\ge\hat{C}N^{-\tau_1+d+\gotd +r+2}\right\}
\end{equation}
where $\hat{C}$ is the constant appearing in Lemma \ref{secondocontrollo}.

\begin{lemma}\label{usofubini}
For all $|j_{0}|\le  6 N$ one has
$ \meas(\calF_{N}(j_{0})) = O ( N^{-d-1}) $. 
\end{lemma}

\prova
By Fubini Theorem we have
$$
\meas(\BBB_{2,N}^{0}(j_{0},\e))=\int_{1/2}^{3/2}{\rm d}\la\,
\meas(\bar B_{2,N}^{0}(j_{0},\e,\la)).
$$
Now, for any $\beta > 0 $, using Lemma \ref{lem.misurabione.nlw} we have
$$
\begin{aligned}
CN^{-\tau_1+d+\nu+1}&\ge\int_{1/2}^{3/2}{\rm d}\la\,
\meas(\bar B_{2,N}^{0}(j_{0},\e,\la))\\
&\ge\beta \meas(\{\la\in[1/2,3/2]\;:\;\meas(\bar B_{2,N}^{0}(j_{0},\e,\la))\ge\beta \})
\end{aligned} 
$$
and for $\beta = \hat{C}N^{-\tau_1+2d +\nu+2}$  we prove the lemma (recall \eqref{FNj0}).
\EP

\begin{lemma}\label{misuralambdabarra} If $ \tau_0 > d+ 3\nu+ 1 $ then 
$ \meas([1/2,3/2]\setminus\ol{\mathcal I})= O(N_{0}^{-1} )$
where $\ol{\mathcal I}$ is defined in \eqref{lambdabarra}.
\end{lemma}

\prova 
Let us write
$$
[1/2,3/2]\setminus\ol{\mathcal I}=\bigcup_{|\ell|,|j|\le N_{0}}
\RR_{l,j},\qquad \RR_{\ell,j}:= \Big\{\la\in\calI\;:\; |(\la\ol{\oo}\cdot \ell)^{2}-|j|^2|\le N_{0}^{-\tau_{0}} 
\Big\}.
$$
Since $j\in\ZZZ^d_*$, then $\RR_{0,j}=\emptyset$ if $N_{0}>1$. For $ \ell \neq 0 $,
 using the Diophantine condition \eqref{dio}, we get  $ \meas(\RR_{\ell,j})\le C N_{0}^{-\tau_{0}+2\nu}, $
so that
$$
\meas([1/2,3/2]\setminus\ol{\mathcal I})\le \sum_{|\ell|,|j|\le N_{0}}\meas(\RR_{\ell,j})\le
C N_{0}^{-\tau_{0}+d+3\nu} = O( N_0^{-1})
$$
because $\tau_{0}-d-3\nu>1$. 
\EP

The measure of the set $ \tilde{\mathcal I} $ in \eqref{lambdatilte} is estimated in \cite{BB2}-Lemma 6.3 
(where $\tilde{\mathcal I}$ is denoted by $\tilde{\mathcal G}$).

\begin{lemma}[]\label{lem.lambdatilde}
If $\g<\min(1/4,\g_{0}/4)$ (where $ \g_0 $ is that in \eqref{diophquad}) then 
$ \meas([1/2, 3/2] \setminus\tilde{\mathcal I})=O(\g)$.
\end{lemma}

To conclude the measure esitimate we note that
by the definition in \eqref{FNj0} 
 for all $\la\not\in\calF_{N}(j_{0})$
one has $\meas(\bar B_{2,N}^{0}(j_{0},\e,\la))<O(N^{-\tau_1+2d+\nu+2})$. Thus for any
$\la\not\in\calF_{N}(j_{0})$, applying Lemma \ref{secondocontrollo} we have
\begin{equation}\nonumber 
\bar B_{N}^{0}(j_{0},\e,\la)\subset\bigcup_{q=1}^{N^{2d+\nu+4}}I_{q} \, , \quad
I_{q} \mbox{ intervals with } \meas(I_q) \le N^{-\tau_1} \, .
\end{equation}
But then, using also Lemma \ref{primapartizione}, we have that (recall \eqref{buoniautovaloriinf}
with $ \gote = 2d+\nu+ 4 $)
\begin{equation}\nonumber 
[1/2, 3/2]  \setminus \bar  \calG_{N}^{0}
\subset\bigcup_{|j_{0}|\le ( \pippa + 5 ) \pippa^{-1}N}\calF_{N}(j_{0}) \, .
\end{equation}
Hence,  using Lemma \ref{usofubini},   
\begin{equation}\nonumber 
{\rm meas} ({\mathcal I} \setminus \bar  \calG_{N}^{0} )\le\sum_{|j_{0}|\le 6N}
\meas(\calF_{N}(j_{0}))\le O(N^{-1}) .
\end{equation}
Moreover by \eqref{cattivinormalta} with  $ \tau_1 > d+ \nu + 2 $ we get
\begin{equation}\label{meas.bad1}
{\rm meas}({\mathcal I} \setminus {\bar \gotG}_{N} ) 
= O( N^{-1} ) ,
\end{equation}
and finally, Lemmas \ref{misuralambdabarra} and \ref{lem.lambdatilde}   with  $ \g = N_0^{-1} $ imply
$$
 {\rm meas}({\mathcal I} \setminus (\ol{\mathcal I}\cap\tilde{\mathcal I})) = O (N_0^{- 1 }) \, . 
 $$
 
 Putting these estimates together and recalling the definition \eqref{Cantorfinale} of $\CCCC_\e$, we have that
\begin{equation}\label{asintotica.nls}
\begin{aligned}
\meas(\mathcal I\setminus\CCCC_{\e})& = \meas\Big(
\bigcup_{n\ge 0}(\bar \calG^{0}_{N_{n}})^{c} \bigcup_{n \ge 0}(\bar \gotG_{N_{n}} )^{c}
\cup\tilde{\mathcal I}^{c}\cup\ol{\mathcal I}^c\Big)\\
&\le \sum_{n\ge0}\meas ({\mathcal I} \setminus \bar \calG^{0}_{N_{n}})  +\sum_{n\ge 0}\meas(  {\mathcal I} \setminus  
\bar \gotG_{N_{n}})
+\meas( {\mathcal I} \setminus (\ol{\mathcal I} \cap \tilde {\mathcal I}))\\
&\stackrel{\eqref{meas.bad1}}{\lesssim}
\sum_{n\ge0}N_n^{- 1}+  N_{0}^{- 1} 
\lesssim N_{0}^{-1}  {\lesssim}  \e^{1/(S+1)}
\end{aligned}
\end{equation}
i.e. $\CCCC_\e$ has asymptotically full measure.
\EP


\end{document}